\documentclass[aap]{imsart}

\usepackage[]{amsmath}
\usepackage[]{amssymb}
\usepackage{xspace}
\usepackage{pstricks}
\usepackage{curves}
\usepackage{amsthm}
\usepackage{natbib}

\usepackage{hyperref}
\usepackage{xr}
\usepackage{listings}
\usepackage{subfigure}
\usepackage{epsfig}
\usepackage[all]{xy}

\newcommand{\freestar}{ \framebox[7pt]{$\star$} }
\newcommand{\boxslash}{\begin{picture}(9,5)  \put(2,0){\framebox(5,5){$\smallsetminus$}} \end{picture}}

\font\BBb = dsrom10

\newcommand{\C}{{\BBb C}}
\newcommand{\Z}{{\BBb Z}}
\newcommand{\N}{{\BBb N}}

\startlocaldefs
\numberwithin{equation}{section}
\theoremstyle{plain}
\newtheorem{teo}{Theorem}[section]
\newtheorem{defin}{Definition}[section]
\newtheorem{lem}{Lemma}[section]
\endlocaldefs
\begin{document}
\bibliographystyle{plain}

\begin{frontmatter}

\title{Multiplicative free Convolution and Information-Plus-Noise Type Matrices}
\runtitle{Multiplicative free Convolution}

\begin{aug}
  \author{\O yvind Ryan\thanksref{t2}\ead[label=e1]{oyvindry@ifi.uio.no}\ead[label=u1,url]{http://www.ifi.uio.no/$\sim$oyvindry}}
  \and
  \author{M{\'e}rouane Debbah\ead[label=e2]{debbah@eurecom.fr}\ead[label=u2,url]{http://www.eurecom.fr/$\sim$debbah/}}

  \thankstext{t2}{Partially sponsored by the project IFANY (INRIA), ACI MALCOM (CNRS) and  the Institute for Mathematical Sciences, National University of Singapore}

  \runauthor{\O yvind Ryan et al.}

  \affiliation{University of Oslo and Institut Eurecom}

  \address{Department of Informatics\\
           Gaustadalleen 23, P.O Box 1080 Blindern\\
           NO-0316 Oslo, Norway,\\
           Mobile Communications Group\\
       2229 Route des Cretes, B.P. 193\\
           06904 Sophia Antipolis Cedex, France,\\
          \printead{e1,e2},\\
          \printead{u1,u2}}
\end{aug}

\begin{abstract}
Free probability and random matrix theory has shown to be a
fruitful combination in many fields of research, such as digital
communications, nuclear physics and mathematical finance. The link
between free probability and eigenvalue distributions of random
matrices will be strengthened further in this paper. It will be
shown how the concept of multiplicative free convolution can be
used to express known results for eigenvalue distributions of a
type of random matrices called Information-Plus-Noise matrices. 
The result is proved in a free probability framework,
and some new results,  useful for problems related to free
probability, are presented in this context. The connection between
free probability and estimators for covariance matrices is also
made through the notion of free deconvolution.
\end{abstract}



\end{frontmatter}

\maketitle

%\date{}

%\lstset{language=c,numbers=left,numberstyle=\tiny,stepnumber=5,numbersep=5pt,stringstyle=\ttfamily,basicstyle=\ttfamily,basewidth=0.45em,showstringspaces=false}

\section{Introduction}
Applications of free probability have been growing rapidly over the last years.
Random matrices and their limit eigenvalue distributions is an area where free probability has proved to be useful~\citep{book:hiaipetz}.
Random matrices are a useful tool for modelling systems,
for instance in digital communications~\citep{paper:telatar99,paper:tsehanly},
nuclear physics~\citep{paper:guhr,book:mehta} and mathematical
finance~\citep{book:bouchaud}.
This paper is a contribution to the random matrix facet of free probability,
in that the connection between certain random matrices and free probability is clarified further.
We will focus on what we call {\em Information-Plus-Noise Type Matrices}, i.e. random matrices on the form
\begin{equation} \label{system1}
  W_n = \frac{1}{N}(R_n + \sigma X_n)(R_n + \sigma X_n)^{\ast},
\end{equation}
where $R_n$ and $X_n$ are independent random matrices of dimension $n\times N$.
These can be thought of as sample covariance matrices of random vectors $r_n + \sigma x_n$, where
$r_n$ can be interpreted as a vector carrying the information in a system,
and $x_n$ additive noise, with $\sigma$ the strength of the noise.
We impose no assumption on independence between samples.
We will use some common restrictions on the noise:
$X_n$ will contain i.i.d. complex entries of unit variance.
$n$ and $N$ will be increased so that
\begin{equation} \label{gcondition}
  \lim_{n\rightarrow\infty}\frac{n}{N} = c.
\end{equation}

In~\citep{paper:doziersilverstein1}, Dozier and Silverstein explain
how the limit eigenvalue distribution $\mu_W$ of the matrix $W_n$
can be found, based on knowledge of the limit eigenvalue
distribution $\mu_{\Gamma}$ of the matrix $\Gamma_n =
\frac{1}{N}R_nR_n^{\ast}$. The result is expressed in terms of a
solution to a function equation (equation (\ref{dozsilv11})). We
will show that there is an equivalent way of expressing this
solution, using the concept of {\em multiplicative free
deconvolution}, denoted by $\boxslash$ 
(multiplicative free convolution, as well as freeness and asymptotic freeness are defined in 
section~\ref{preliminaries}). The following is the main result of
the paper:
\begin{teo} \label{teo1d}
  Assume that the entries $X_{ij}^n$ of $X_n$ are Gaussian, independent and identically distributed with expectation 0 and variance 1.
  Assume also that the empirical eigenvalue distribution of $\Gamma_n=\frac{1}{N}R_nR_n^{\ast}$ converges in distribution
  almost surely to a compactly supported probability measure $\mu_{\Gamma}$.
  Then we have that the empirical eigenvalue distribution of $W_n$ also converges in distribution almost surely
  to a compactly supported probability measure $\mu_W$ uniquely identified by
  \begin{equation} \label{teo1dequation}
    \mu_W \boxslash \mu_c = (\mu_{\Gamma} \boxslash \mu_c) \boxplus \mu_{\sigma^2 I}.
  \end{equation}
\end{teo}
Some remarks are needed to explain theorem~\ref{teo1d}.
By the {\em empirical eigenvalue distribution}
of an $n\times n$ random matrix $X$ we mean the (random) atomic measure
\[
  \frac{1}{n} \left( \delta(\lambda_1(X)) + \cdots + \delta(\lambda_n(X)) \right),
\]
where $\lambda_1(X),...,\lambda_n(X)$ are the (random) eigenvalues of $X$.
That $\mu_n$ converges in distribution to $\mu$ means that the moments of $\mu_n$ converge to the moments of $\mu$.
Theorem~\ref{teo1d} requires compactly supported measures, and these have moments of all orders.

The conditions in theorem~\ref{teo1d} are somewhat stronger than those in~\citep{paper:doziersilverstein1}
due to the restriction to measures with compact support.
Contrary to~\citep{paper:doziersilverstein1}, we also restrict to noise-matrices with Gaussian entries.

Theorem~\ref{teo1d} yields a short expression for $\mu_W$,
removing the need for solving the equation
in~\citep{paper:doziersilverstein1} directly. 
It essentially says that the connection between $\mu_W$ and $\mu_{\Gamma}$ 
can be expressed compactly in the {\em deconvolved domain}, where the connection can be viewed 
as a shift of the spectrum with the noise variance $\sigma^2$.
The proof of theorem~\ref{teo1d} is based on
methods from free probability, with some new  results established
along the way. Some of these deserve extra attention, in
particular theorem~\ref{prop2}. This can be thought of as a
version of theorem~\ref{teo1d} where $\mu_W$ and $\mu_{\Gamma}$ are interpreted as
distributions of free random variables.

Theorems~\ref{23101} and~\ref{23102} also deserve some extra attention.
These address asymptotic freeness almost everywhere~\citep{book:hiaipetz} for two random matrices where
\begin{enumerate}
  \item both converge in distribution almost everywhere to compactly supported limits, and
  \item one of the random matrices are standard unitary (theorem~\ref{23101}) or Gaussian (theorem~\ref{23102}).
\end{enumerate}
These results expand known results from~\citep{book:hiaipetz} for
asymptotic freeness. The proofs of theorems~\ref{23101}
and~\ref{23102} use random matrix approximations with
deterministic matrices. Asymptotic freeness of Gaussian/standard
unitary random matrices and uniformly norm-bounded deterministic
matrices are well-known (lemma 4.3.2 in~\citep{book:hiaipetz}).
Unfortunately, norm-bounded deterministic matrices are not able to
approximate the random matrices under consideration. We solve the
problem by generalizing to matrices satisfying uniform $\|\cdot \|_p$-norm
bounds instead, where $\|\cdot \|_p$ dentotes the Schatten $p$-norm (with respect to $tr_n$), 
defined for $p\geq 1$ by $\| A \|_p = tr_n( | A |^p)^{\frac{1}{p}}$ ($A\in M_n(\C)$): 
We prove that matrices satisfying such bounds can be used to
approximate our random matrices, and that they also give
asymptotic freeness as in lemma 4.3.2 in~\citep{book:hiaipetz}
(theorem~\ref{432}).

Theorem~\ref{teo1d} is actually proved by combining theorems~\ref{23102} and~\ref{prop2}
through another approximation argument (see theorem~\ref{prop2b}).
While~\citep{paper:doziersilverstein1} restricts to the distribution of $\frac{1}{N}(R_n+\sigma X_n)(R_n+\sigma X_n)^{\ast}$,
we show more in that any mixed moments of $\frac{1}{N}R_nR_n^{\ast}$ and $\frac{1}{N}X_nX_n^{\ast}$
are obtained through our asymptotic freeness results.

Recent works~\citep{eurecom:freedeconvinftheory,fdsigpro} show that
multiplicative free convolution also admits an efficient
implementation in terms of the moments of the operand measures. 
The basic results on free probability we need for
this are proved in this paper (theorems~\ref{betterform0}
and~\ref{computableform0}). A consequence is that existing
computational frameworks can be used in obtaining $\mu_{\Gamma}$ and $\mu_W$.
In~\citep{eurecom:freedeconvinftheory}, $\mu_{\Gamma}$ and $\mu_W$
are illustrated in terms of signal processing applications, and
simulations are run using a computational framework building on
theorems~\ref{betterform0} and~\ref{computableform0}. A useful consequence of the link with free
probability is that that the "inverse problem" (i.e. that of
finding $\mu_{\Gamma}$ from $\mu_W$) can be solved within the same
framework, since the framework embraces 
convolution as well as deconvolution. 

The eigenvalue distribution of $\Gamma_n$ provides us with
possibilities for estimating the covariance matrices of the system
through the so-called {\em $G$-estimators} ~\citep{book:girkostat}.
These will be reviewed, and it will be shown how multiplicative
free convolution can be used to rewrite such estimators to a very
simple form. It will be apparent from this that the
$G^2$-estimator actually can be viewed as a step in expressing
$\mu_W$ from $\mu_{\Gamma}$.

While the results mentioned here are hard to prove, some of them
should should come as no surprise. For
instance,~\citep{paper:raoedelman} has already made the connection
between Information-Plus-Noise type matrices and multiplicative
free convolution. This paper also indicates that some of the
mentioned results are already known, by saying that random
matrices with Haar-distributed eigenvectors are asymptotically
free from any random matrices independent from them. However, the
generality in which this should hold is not indicated.
Also,~\citep{paper:raoedelman} considers only Gaussian matrices,
and the connection with already existing estimators of covariance
matrices was not made.

This paper is organized as follows. Section~\ref{preliminaries}
contains notation and preliminaries for various free probability
tools, like free transforms and combinatorial aspects. The
mentioned implementation of free convolution builds on the
combinatorial expression of freeness, and the results needed on
this are explained in section~\ref{frameworkimpl}. 
The proof for theorem~\ref{teo1d} is presented in section~\ref{reproof}. 
A sketch of the proof is first given, 
followed by the proofs for theorems~\ref{432},~\ref{23101}~\ref{23102} and~\ref{prop2}.
Section~\ref{sectionsystem1} first states the results we need from~\citep{paper:doziersilverstein1}, 
and sketch the proof for the equivalence of these and theorem~\ref{teo1d}. 
This sketch is then followed by the rest of the details. 
The various transforms used in free probability (section~\ref{preliminaries})
are used in this direction. 
In section~\ref{ganalysis} we state the principles of $G$-analysis
and the expression for the $G^2$-estimator. We also prove
the theorem which expresses the $G^2$-estimator in terms of free
probability.
\section{Notation and preliminaries} \label{preliminaries}
In the following, uppercase symbols will be used for matrices, and 
$(.)^{\ast}$ will denote hermitian transpose. $I_n$ will represent the
identity matrix of order $n$. We will focus here on certain
noncommutative probability spaces. A noncommutative probability
space is a pair $(A,\phi)$ where $A$ is a unital $\ast$-algebra
and $\phi$ is a normalized (i.e. $\phi(I)=1$) linear functional on
$A$. The elements of $A$ are called random variables. The
probability spaces we will encounter are mostly $(M_n(\C),tr_n)$,
i.e. $n\times n$-matrices equipped with the normalized trace. Any
matrix can be associated with a probability measure through it's
eigenvalue distribution. We will mostly be concerned with
probability measures with compact support.

\begin{defin} \label{freedef}
  A family of unital $\ast$-subalgebras
  $(A_i)_{i\in I}$ will be called a free family if
\begin{equation}
       \left\lbrace  \begin{matrix}
           a_j\in A_{i_j} \\
           i_1\neq i_2,i_2\neq i_3,\cdots ,i_{n-1}\neq i_n \\
           \phi(a_1)=\phi(a_2)=\cdots =\phi(a_n)=0 \end{matrix}
       \right\rbrace
  \Rightarrow \phi(a_1\cdots a_n)=0.
\end{equation}
\end{defin}
(\ref{freedef}) enables us to calculate the mixed moments of $a_1$ and $a_2$ when they are free.
In particular, the moments of $a_1+a_2$ and $a_1a_2$ can be calculated.
This gives us two new probability measures, which depend on the
probability measures of $a_1$, $a_2$ only (i.e. not on their realizations).
Therefore we can define two operations on the set of probability measures:
Additive free convolution
\begin{equation} \label{addconvdef}
  \mu_1\boxplus\mu_2
\end{equation}
for the sum of free random variables, and multiplicative free convolution
\begin{equation} \label{multconvdef}
  \mu_1\boxtimes\mu_2
\end{equation}
for the product of free random variables.

Let $F^{\mu_A}$ denote the empirical distribution function (e.d.f.) of the eigenvalues of $A$
(so that $F^{\mu_A}(x)$ is the proportion of eigenvalues of $A$ which are $\leq x$).
When we have a series of e.d.f.'s $F^{\mu_{A_n}}$,
we will use the notation
\[
  F^{\mu_{A_n}} \stackrel{\cal D}{\rightarrow} F^{\mu}
\]
for weak convergence, where $F^{\mu}$ is the cumulative distribution function of the measure $\mu$.
We will also write a.s. as shorthand notation for almost sure convergence.

Some random matrices and limit distributions occur naturally in many contexts.
If the entries of the $n\times N$ (with $\lim_{n\rightarrow\infty} \frac{n}{N} = c$) random matrices $W_n$ have zero mean and unit variance,
the empirical eigenvalue distribution of
$\frac{1}{N} W_nW_n^{\ast}$ converges almost surely to the so-called Mar\u{c}henko Pastur law $\mu_c$
(\citep{book:tulinoverdu} page 9).
These are also called the free Poisson distributions, and are characterized by the density
\begin{equation} \label{mpdensity}
  f^{\mu_c}(x) = (1-\frac{1}{c})^+ \delta(x) + \frac{\sqrt{(x-a)^+(b-x)^+}}{2\pi cx},
\end{equation}
where $(z)^+ =\mbox{max}(0,z)$, $a=(1-\sqrt{c})^+$ and $a=(1+\sqrt{c})^+$. 
Similar notation to the Mar\u{c}henko Pastur law is used for the distribution $\mu_a$ of a random variable $a$. 
We avoid confusion by never using $c$ to denote random variables. 
$\mu_1$ will always mean the Mar\u{c}henko Pastur law with parameter one. 
Mar\u{c}henko Pastur laws are some of the most basic random matrix building blocks, as they appear as limits for large random matrices in many contexts.
This paper will demonstrate that this is indeed the case for the type of systems we consider also.

We will not use the characterization of the Mar\u{c}henko Pastur law as in (\ref{mpdensity}) directly.
Rather we will work with equivalent expressions of it through the transforms defined in this section.
The transforms we define will only be applied for probability measures with support contained on the positive real line.

{\em The Stieltjes transform} (\citep{book:tulinoverdu} page 38) of a probability measure $\mu$ is the analytic function on $\C^+ = \{ z\in C : \Im z > 0 \}$ defined by
\begin{equation}
  m_{\mu}(z) = \int_{-\infty}^{\infty} \frac{1}{\lambda-z} dF^{\mu}(\lambda).
\end{equation}
A convenient inversion formula for the Stieltjes transform also exists, so that $m_{\mu}$ uniquely identifies $\mu$.
If $\mu$ is assumed to have nonnegative support, $m_{\mu}$ can be analytically continued to the negative part of the real line.
If $\mu = \mu_X$ for a non-negative random variable $X$,
$m_{\mu}$ is strictly monotone on the negative real line,
taking values in the interval $[ 0,E\left( \frac{1}{X}\right) ]$.
We will use the fact that if we know $m_{\mu}(z)$ in an interval $(-z,0)$ for $z<0$, we also know $m_{\mu}$ for all other values of $z$,
and hence we also know $\mu$ (use the Stieltjes inversion formula).

{\em The $\eta$-transform} (\citep{book:tulinoverdu} page 40) is defined for measures $\mu$ with support on the positive real line,
and for nonnegative real numbers by
\begin{equation}
  \eta_{\mu}(z) = \int_{-\infty}^{\infty} \frac{1}{1+z\lambda} dF^{\mu}(\lambda).
\end{equation}
$\eta(z)$ is a strictly monotonically decreasing function. As such
it simplifies many derivations and statements of results. The
inverse is tightly connected to the $S$-transform (see below).
It's connection with the Stieltjes transform is
\begin{equation}
  \eta_{\mu}(z) = \frac{m_{\mu}(-\frac{1}{z})}{z}, \mbox{   } m_{\mu}(z) = -\frac{\eta_{\mu}(-\frac{1}{z})}{z}.
\end{equation}
Therefore $\eta_{\mu}(z)$ uniquely identifies $m_{\mu}(z)$, since $m_{\mu}(z)$ for real, negative $z$ can be continued analytically to $\C^+$.
We will use the fact that if we know $\eta_{\mu}(z)$ in an interval $(0,z)$ for $z>0$, we also know $\mu$.

{\em The $R$-transform} (\citep{book:tulinoverdu} page 48) has domain of definition $\C^+$ and can be defined in terms of the Stieltjes transform as
\begin{equation}
  {\cal R}_{\mu}(z) = m_{\mu}^{-1}(-z) - \frac{1}{z}.
\end{equation}
The importance of the $R$-transform comes from it's additive property for the distribution of the sum of free random variables $A_1$ and $A_2$,
\begin{equation} \label{radditivity}
  {\cal R}_{\mu_{a_1+a_2}}(z) = {\cal R}_{\mu_{a_1}}(z) + {\cal R}_{\mu_{a_2}}(z).
\end{equation}
Slightly different versions of the $R$-transform are encountered in the litterature. The one above is from~\citep{book:tulinoverdu}.
In connection with free combinatorics, another definition is used, namely $R_{\mu}(z) = z{\cal R}_{\mu}(z)$.
Of course, $R_{\mu}(z)$ also satisfies (\ref{radditivity}).

{\em The $S$-transform} (\citep{book:tulinoverdu} page 50) is defined on $(-1,0)$. It can be defined in terms of the $\eta$-transform by
\begin{equation}
  S_{\mu}(z) = -\frac{z+1}{z}\eta_{\mu}^{-1}(z+1).
\end{equation}
The Mar\u{c}henko Pastur law (\ref{mpdensity}) can be shown to have S-transform $S_{\mu_c}(z) = \frac{1}{1+cz}$ (\citep{book:tulinoverdu} page 51).
The importance of the $S$-transform comes from it's multiplicative property for the distribution of the product of free random variables $a_1$ and $a_2$:
\begin{equation}
  S_{\mu_{a_1 a_2}}(z) = S_{\mu_{a_1}}(z) S_{\mu_{a_2}}(z).
\end{equation}
If the values of $\eta_{\mu}(z)$ or $S_{\mu}(z)$ are known in an interval, one also knows $\mu$.

\begin{figure}
  \begin{center}
    \epsfig{figure=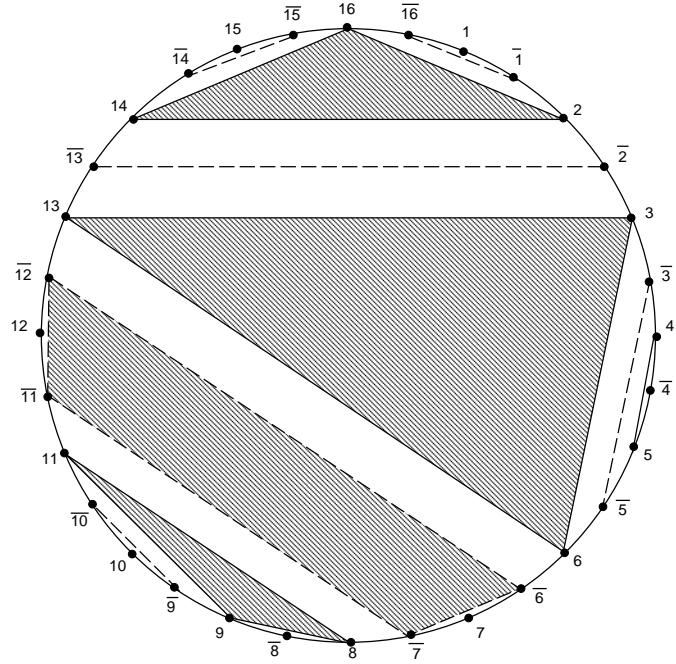,width=.7\columnwidth}
  \end{center}
  \caption{The circular representation of a partition of $\{ 1,...,16\}$.}\label{fig:big}
\end{figure}

Freeness, additive and multiplicative free convolution have a
combinatorial description involving these transforms which we will
 use for in some of our proofs. These combinatorial descriptions
build on the concept of {\em noncrossing partitions}:
\begin{defin}
  A partition $\pi$ is called noncrossing if whenever we have $i<j<k<l$ with $i\sim k$, $j\sim l$ ($\sim$ meaning belonging to the same block),
  we also have $i\sim j\sim k\sim l$ (i.e. $i,j,k,l$ are all in the same block).
  The set of noncrossing partitions of $\{ 1,,,.,n\}$ is denoted $NC(n)$.
\end{defin}
$NC(n)$ becomes a lattice under the refinement order of partitions.
An ingredient we need in making the connections between freeness and the noncrossing partitions is the complementation map of Kreweras,
which is a lattice anti-isomorphism of $NC(n)$.
To define this we need the circular representation of  a partition: We mark $n$ equidistant points $1,...,n$ (numbered clockwise) on the circle,
and form the convex hull of points lying in the same block of the partition.
This gives us a number of convex sets $H_i$, equally many as there are blocks in the partition, which do not intersect if and only if the partition is noncrossing.
Put names $\bar{1},...,\bar{n}$ on the midpoints of the $1,...,n$ (so that $\bar{i}$ is the midpoint of the segment from $i$ to $i+1$).
The complement of the set $\cup_i H_i$ is again a union of disjoint convex sets $\tilde{H}_i$.
All this is demonstrated in figure~\ref{fig:big},
where the $\tilde{H}_i$ are the scrambled areas with dashed borders,
the $H_i$ are the scrambled areas with non-dashed borders.
We will refer to this figure heavily during the proof of theorem~\ref{prop2}.
We can now define the Kreweras complementation map:
\begin{defin} \label{krewerasdef}
  The Kreweras complement of $\pi$, denoted $K(\pi)$, is the partition on $\{ \bar{1},...,\bar{n}\}$ determined by
  \[
    i\sim j \mbox{ in } K(\pi) \iff \bar{i},\bar{j} \mbox{ belong to the same convex set } ~\tilde{H}_k.
  \]
\end{defin}
The connection between the $R$-transform and noncrossing partitions comes through the {\em moment-cumulant formula},
which relates the moments and the $R$-transform coefficients (also called cumulants) for the distribution of a random variable.
\begin{lem}
  Write the $R$-transform as a power series, $R_{\mu_a}(z) = \sum_n \alpha_n z^n$. Then
  \begin{equation}
   \phi(a^n) = \sum_{\pi = \{ B_1,\cdots ,B_k \}\in NC(n)} \prod_{i=1}^k \alpha_{|B_i|},
  \end{equation}
\end{lem}
This can be used as an alternative definition of the
$R$-transform. We also need to define the multidimensional
$R$-transform for the joint distribution of a sequence of random variables. Denote by $C\langle
z_1,...,z_m\rangle$ the space of complex power series in $m$
noncommuting variables $z_i$ with vanishing constant term. These
can be written in the form
\[
  \sum_{k\geq 1} \sum_{i_1,...,i_k} a_{i_1,...,i_k} z_{i_1} \cdots z_{i_k}.
\]
In referring to the coefficients of a power series $f$ on this form we will write
\[
  [ coef(i_1,...,i_k) ](f) = a_{i_1,...,i_k},
\]
and if $\pi = \{ B_1,...,B_m\}$,
\[
\begin{array}{rcl}
  [ coef(i_1,...,i_k)|B_i  ](f) &=& a_{{(i_j)}_{j\in B_i}}\\
  \left[ coef(i_1,...,i_m);\pi \right] (f) &=& \prod_i [ coef(i_1,...,i_k)|B_i ](f).
\end{array}
\]
For power series in one variable, the coefficients will also be
written in the form $[ coef_k ](f)$. For $n$ random variables
$a_1,...,a_n$ we define their joint moment series as the power
series $M_{\mu_{a_1,...,a_n}}\in C\langle z_1,...,z_n\rangle$ such
that
\[
  M_{\mu_{a_1,...,a_n}}(z_1,...z_k) = \sum_{m\geq 1} \sum_{i_1,...,i_m} \phi(a_{i_1}\cdots a_{i_m}) z_{i_1}\cdots z_{i_m},
\]
and we define their joint $R$-series as the unique power series $R_{\mu_{a_1,...,a_n}}\in C\langle z_1,...,z_n\rangle$ such that
\begin{equation} \label{jointrseriesdef}
  \phi(a_{i_1}\cdots a_{i_m}) = \sum_{\pi\in NC(m)} [ coef(i_1,...,i_m);\pi ](R_{\mu_{a_1,...,a_n}}).
\end{equation}
The result we will use connecting the joint $R$-series and freeness is the following:
\begin{lem} \label{nomixedterms}
  $( \{ a_1,...,a_n\} , \{ b_1,...,b_m\} )$ is a free family if and only if
  \[
    R_{\mu_{a_1,...,a_n,b_1,...,b_m}}(z_1,...,z_{n+m})
  \]
  \[
    = R_{\mu_{a_1,...,a_n}}(z_1,...,z_n) + R_{\mu_{b_1,...,b_m}}(z_{n+1},...,z_{n+m}).
  \]
\end{lem}
This lemma is often summarized by saying that the joint $R$-series
of free random variables has no {\em mixed} terms. A  special form
of (\ref{jointrseriesdef}) and lemma~\ref{nomixedterms} we will
use is the following: If $(a_1,a_2)$ is a free family, and a mixed
term $a_{i_1}\cdots a_{i_m}$ is given, form the partition with two
blocks $\sigma = \{ \sigma_1 , \sigma_2 \}$, where $\sigma_k = \{
j | a_{i_j} = k \}$. Then
\begin{equation} \label{usefulsum}
  \phi(a_{i_1}\cdots a_{i_m}) = \sum_{\pi\le\sigma\in NC(m)} [ coef(i_1,...,i_m);\pi ](R_{\mu_{a_1,a_2}}).
\end{equation}
Our combinatorial connection with multiplicative free convolution can be made complete with the help of the following definition~\citep{nica1},~\citep{ns2}:
\begin{defin}
  Given two power series $f$ and $g$, their {\em boxed convolution} $f\freestar g$ is defined by
  \[
    [ coef(i_1,...,i_m) ](f\freestar g)
  \]
  \begin{equation} \label{boxedconvdef}
    = \sum_{\pi\in NC(m)} [ coef(i_1,...,i_m);\pi ](f) [ coef(i_1,...,i_m);K(\pi) ](g).
  \end{equation}
\end{defin}
Boxed convolution is commutative only on power series in one variable~\citep{book:comblect}.
It satisfies the associative law, but not the distributive law. It does not satisfy linearity properties w.r.t. scalar multiplication.
However, the following holds and will be useful to us:
\[
  [ coef_n ] \left( (cf)\freestar (cg) \right) = [ coef_n ] \left( c^{n+1}(f\freestar g) \right)
\]
and
\[
  [ coef_n ] \left( f\freestar (cId) \right) = [ coef_n ] \left( c^n f \right).
\]
Here we used the shorthand notation $c^n f$ for the power series defined by $[ coef_n ] (c^nf) = c^n [coef_n ](f)$.
The first statement is easily proved using the fact that $|\pi |+|K(\pi)| = n+1$ for any $\pi\in NC(n)$~\citep{book:comblect}.
The second statement is trivial.
The following result holds for multiplicative free convolution~\citep{book:comblect}:
\begin{lem}
  If $( \{ a_1,...a_n\} , \{ b_1,...,b_n\})$ is a free family, then
  \begin{equation} \label{boxedconv1}
    R_{\mu_{a_1b_1,...,a_nb_n}} = R_{\mu_{a_1,...,a_n}} \freestar R_{\mu_{b_1,...,b_n}}
  \end{equation}
\end{lem}
One can also define additive and multiplicative free deconvolution in most cases,
i.e. finding $\mu_2$ in (\ref{multconvdef}) when $\mu_1\boxtimes\mu_2$ are known.
\begin{defin}
  Given probability measures $\mu$ and $\mu_2$.
  When there is a unique probability measure $\mu_1$ such that $\mu = \mu_1 \boxtimes \mu_2$,
  we will denote $\mu_1 = \mu \boxslash \mu_2$.
  We say that $\mu_1$ is the multiplicative free deconvolution of $\mu$ with $\mu_2$.
\end{defin}
We can define addtive free deconvolution similarly. 
Note that free deconvolution is defined only for a subset of all probability measures, 
since measures exist which can't be expressed on the forms $\mu_1\boxplus\mu_2$ or $\mu_1\boxtimes\mu_2$. 
Deconvolution can, however, also be viewed as a formal operation on a sequence of moments. 
Viewed as such, multiplicative free deconvolution is well-defined when we have non-vanishing first moments. 
This can be seen from the combinatorial description of multiplicative free
convolution (\ref{boxedconvdef}). In light of (\ref{boxedconv1}),
it is obvious from (\ref{boxedconvdef}) that the cumulants in
$R_{\mu_2}$ can be calculated recursively from those of $R_{\mu_1}$
and $R_{\mu_1\boxtimes\mu_2}$, when the first coefficient of
$R_{\mu_1}$ (which equals the first moment) is known. 
Since the main theorem relates to the moments of the involved measures 
(it is a statement on convergence in distribution), we will in the following view deconvolution in terms of the moments only. 

A form of (\ref{boxedconv1}) which will be useful to us is for the case $n=1$. If we write $\mu_a = (\mu_a\boxslash\mu_c)\boxtimes\mu_c$, we get
\begin{equation} \label{boxedconvusefulform}
  R_{\mu_a\boxslash\mu_c} = R_{\mu_a} \freestar R_{\mu_c}^{-1}.
\end{equation}
The facts we will use concerning boxed convolution are the following,
relating moment series, $R$-series of general random variables (in
particular projections and free Poisson random variables), the
$Zeta$ series, which is defined as
\[
  Zeta(z_1,...,z_n) = \sum_k \sum_{i_1,...,i_k} z_{i_1} \cdots z_{i_k},
\]
the $Moeb$ series (which is the inverse of $Zeta$ under composition with $\freestar$)
and $Id$ (which is the unit under composition with $\freestar$):
\begin{equation} \label{cconnection}
  \begin{array}{lclclcl}
    M_{\mu}   &=& R_{\mu} \freestar Zeta  &\mbox{ and }& R_{\mu}   &=& M_{\mu} \freestar Moeb\\
    M_{\mu_p} &=& c Zeta                  &\mbox{ and }& R_{\mu_c} &=& c^{n-1} Zeta.
  \end{array}
\end{equation}
Here $p$ is a projection with $\phi(p)=c$. Our definition of
$\mu_c$ differs from that of~\citep{book:hiaipetz}, for purposes of
compatibility
with~\citep{paper:doziersilverstein1}~\citep{paper:raoedelman}.
Consequently, the expressions for the $R$-transforms are
different. In the terminology of~\citep{book:hiaipetz}, the
$R$-series would be $c Zeta$. The following definition~\citep{ns3}
will also be in use:
\begin{defin}
  A pair $(a,b)$ of noncommutative random variables is called an {\em $R$-diagonal pair} if it's $R$-series is of the form
  \begin{equation}
    R_{\mu_{a,b}}(z_1,z_2) = \sum_{n=1}^{\infty} \alpha_n \left(  (z_1 z_2)^n + (z_2 z_1)^n \right).
  \end{equation}
  An element $a$ will be said to be an {\em $R$-diagonal element} if $(a,a^{\ast})$ is an $R$-diagonal pair.
  The one-variable series $\sum_{n=1}^{\infty} \alpha_n z^n$ will be called the determining series of the $R$-diagonal pair $(a,b)$.
\end{defin}
We will  use  the fact that if $a$ is an $R$-diagonal element,
it's determining series can be written as $R_{\mu_{aa^{\ast}}}
\freestar Moeb$~\citep{book:hiaipetz}. Two important $R$-diagonal
elements are 
\begin{enumerate}
  \item the {\em Haar unitary}, which can be defined as a
    unitary $u$ satisfying $\phi(u^n) = 0$ for all $n\in Z\neq 0$, and 
  \item the {\em circular element}, which can be defined as an element $s$ 
    whose $\ast$-distribution $\mu_{s,s^{\ast}}$ 
    satisfies $R_{\mu_{s,s^{\ast}}}(z_1,z_2) = z_1z_2+z_2z_1$. 
\end{enumerate}
The concept of $R$-diagonally was in fact invented
in search of a common approach for Haar unitaries and circular
elements~\citep{ns3}. Haar unitaries are very important in
asymptotic random matrix results. In $W^{\ast}$-probability
spaces, when the isometric part of an $R$-diagonal element has
kernel equal to zero, the isometric part is actually a Haar
unitary.

\subsection{Implementation of free convolution} \label{frameworkimpl}
While free convolution has an abstract definition,
the combinatorial description given in this section can actually be used to obtain an efficient implementation.
In many practical cases, free convolution with $\mu_c$ is what we are interested in. Such free convolution is simplified through the following result.
\begin{teo} \label{betterform0}
\begin{equation} \label{betterform}
  (c M_{\mu}) \freestar Zeta = c\left( M_{\mu \boxtimes \mu_c} \right).
\end{equation}
\end{teo}

\begin{proof}
To see this, start by combining (\ref{boxedconv1}) with (\ref{cconnection}) to get
\[
  R_{\mu \boxtimes \mu_c} = R_{\mu} \freestar R_{\mu_c} = R_{\mu} \freestar (c^{m-1} Zeta).
\]
After convolving both sides with $Zeta$, we get
\begin{equation} \label{betterformstart}
  M_{\mu \boxtimes \mu_c} = M_{\mu} \freestar (c^{m-1} Zeta).
\end{equation}
To prove (\ref{betterform}), rewrite the left hand side as
\[
  \sum_{\pi\in NC(m)}   c^{|\pi |} \left[ coef_m;\pi \right] M_{\mu}.
\]
Since $|\pi | + |K(\pi)| = m+1$, this equals
\[
\begin{array}{ll}
   & c\sum_{\pi\in NC(m)}   \left[ coef_m;\pi \right] (M_{\mu}) c^{m-|K(\pi)|}\\
  =& c\sum_{\pi\in NC(m)}   \left[ coef_m;\pi \right] (M_{\mu}) c^{m-|K(\pi)|} \left[ coef_m;K(\pi) \right] (Zeta)\\
  =& c\sum_{\pi\in NC(m)}   \left[ coef_m;\pi \right] (M_{\mu}) c^m \left[ coef_m;K(\pi) \right] (c^{-1} Zeta)\\
  =& c\sum_{\pi\in NC(m)}   \left[ coef_m;\pi \right] (M_{\mu}) \left[ coef_m;K(\pi) \right] (c^{m-1} Zeta)\\
  =& c\left( M_{\mu} \freestar (c^{m-1} Zeta) \right),
\end{array}
\]
substituting (\ref{betterformstart}) proves the claim.
\end{proof}

In summary, if we need to compute the moments of $\mu \boxtimes \mu_c$,
one can first compute the moment series $c M_{\mu}$,
then use this to compute the left hand side of (\ref{betterform}).
According to (\ref{betterformstart}) and (\ref{betterform}),
the moment series of $\mu \boxtimes \mu_c$ can then be computed from this with an additional scaling with $\frac{1}{c}$.

In other words, convolving with $\mu_c$ is equivalent to convolving with $\mu_1$
(with additional scalings of power series taken into account), since $R_{\mu_1} = Zeta$.
It turns out that boxed convolution with $Zeta$ is easy to compute, as the following result shows.
The result is stated in terms of the moment-cumulant formula, since the relation between cumulants and moments are given by
boxed convolution with $Zeta$.

\begin{teo} \label{computableform0}
\begin{equation} \label{computableform}
  [ coef_m ] (M_{\mu}) = \sum_{k=1}^m [ coef_k ] (R_{\mu}) [ coef_{m-k} ] (1+M_{\mu})^k.
\end{equation}
\end{teo}

\begin{proof}
For each $\pi\in NC(m)$,
fix the block $B_1 = \{ b_{11},...,b_{1k} \}$ in $\pi$ containing $1$,
and let $NC(m,B_1)$ be the set of all noncrossing partitions which contain $B_1$ as a block.
Rewrite the definition of boxed convolution (\ref{boxedconvdef}) to
\begin{equation} \label{neweq1}
\begin{array}{lll}
  [ coef_m ] (M_{\mu}) &=& \sum_{B_1} \sum_{\pi\in NC(m,B_1)} [ coef_m;\pi ](R_{\mu}) [ coef_m;K(\pi) ](Zeta)\\
                       &=& \sum_{B_1} \sum_{\pi\in NC(m,B_1)} [ coef_m;\pi ](R_{\mu}).
\end{array}
\end{equation}
Blocks in $\pi\in NC(m,B_1)$ other than $B_1$ must be entirely contained in one of
$\{ b_{11} + 1,...,b_{12} - 1 \} ,..., \{ b_{1k} + 1,...,b_{11} - 1 \}$.
This means that the inner summand in (\ref{neweq1}) can be rewritten to
\begin{equation} \label{lastsum}
  [ coef_k ](R_{\mu})
  \prod_{i=1}^k \left( \sum_{\pi\in NC( b_{1(i+1)} - b_{1i} - 1)} [ coef_m;\pi ](R_{\mu}) \right).
\end{equation}
From the moment-cumulant formula it is seen that each sum here is simply a moment, so we can rewrite to
\[
  [ coef_k ](R_{\mu})
  \prod_{i=1}^k [ coef_{b_{1(i+1)} - b_{1i} - 1} ] (1+M_{\mu}),
\]
where the summand $1$ in $1+M_{\mu}$ accounts for elements $i$ in (\ref{lastsum}) with $b_{1(i+1)} = b_{1i} + 1$ 
(i.e. consecutive elements in a block). 
All in all, (\ref{neweq1}) can be rewritten to
\begin{equation} \label{from1}
  \sum_k \sum_{ \stackrel{B_1}{|B_1| = k} } [ coef_k ](R_{\mu}) \prod_{i=1}^k \left( [ coef_{b_{1(i+1)} - b_{1i} - 1} ] (1+M_{\mu}) \right)
\end{equation}
Write $a_i = b_{1(i+1)} - b_{1i} - 1$, and note that
\[
  \sum_{i=1}^k a_i = \sum_{i=1}^k \left( b_{1(i+1)} - b_{1i} - 1 \right) = m-k.
\]
The $a_i$ are in one-to-one correspondence with all candidates for $B_1$, so that we can rewrite (\ref{from1}) to
\[
  \sum_k  [ coef_k ](R_{\mu}) \sum_{ \stackrel{a_1,...,a_k}{\sum a_i = m-k} } \prod_{i=1}^k \left( [ coef_{a_i} ] (1+M_{\mu}) \right).
\]
The inner sum here is easily recognized as coefficient $m-k$ in the power series $(1+M_{\mu})^k$
(one factor for each $a_i$). Putting things together we get (\ref{computableform}).
\end{proof}

In (\ref{computableform}) we see that there is no reference to noncrossing partitions.
(\ref{computableform}) can be used easily in calculating moments recursively from cumulants.
The coefficients in the power series $(1+M_{\mu})^k$ can be computed in terms of $k$-fold (classical) convolution.
This is done in~\citep{eurecom:freedeconvinftheory}, where many multiplicative free convolutions are computed
based on (\ref{computableform}). The actual implementation of (\ref{computableform}) used in~\citep{eurecom:freedeconvinftheory} is
contained in~\citep{eurecom:freeimpl}.

Free convolution as introduced here is just defined for compactly supported probability measures.
\section{Proof of theorem~\ref{teo1d}} \label{reproof}
In what follows we first sketch the proof of theorem~\ref{teo1d}.
After this follows proofs for theorems needed in the proof.

First we prove the following variant of lemma 4.3.2 in~\citep{book:hiaipetz}, 
which can be used together with the Borel-Cantelli lemma to prove almost sure convergence.
It is slightly more general in the sense that boundedness in the operator norm $\|\cdot \|$ is not assumed, 
only boundedness in $\| \cdot \|_p$, for $p\geq 1$.
This weaker boundedness assumption is needed since uniformly norm-bounded matrices are not sufficient 
to approximate all compactly supported probability measures almost surely. 
Recall that an $n\times n$ unitary random matrix is called {\em standard unitary} if it's distribution equals 
the Haar probability measure on ${\cal U}(n)$. 
\begin{teo} \label{432}
  Let $U(s,n)_{s\in S}$ be an independent family of $n\times n$ standard unitary random matrices.
  Let $s_1,...,s_l\in S$, $m_1,...,m_l\in \Z\setminus\{ 0\}$, and let $R_p\geq 0$, $p\geq 1$ be constants.
  Then
  \begin{equation} \label{432expression}
    E \left( | tr_n\left( U(s_1,n)^{m_1} D_1(n) U(s_2,n)^{m_2} D_2(n) \cdots U(s_l,n)^{m_l} D_l(n) \right) |^2 \right)
  \end{equation}
  is $O(n^{-2})$ as $n\rightarrow\infty$ uniformly for the choice of any $D_r(n)\in M_n(\C)$ ($1\leq r\leq l$) such that for $1\leq r\leq l$ either
  \[
    tr_n(D_r(n)) = 0 \mbox{ and }  \| D_r(n) \|_p \leq R_p \mbox{ } (n\in\N)
  \]
  or
  \[
    D_r(n) = I_n \mbox{ } (n\in\N) \mbox { and } s_r\neq s_{r+1} (\mbox{ with } s_{l+1} = s_1).
  \]
  Also, for a given $l$, there exists a $p_l$ such that the same statement holds as long 
  as the $\| \cdot \|_p$-norm bounds are satisfied for $p\leq p_l$ only. 
\end{teo}
The proof is in section~\ref{proof432}. It somewhat simplifies the proof of lemma 4.3.2 in~\citep{book:hiaipetz},
and can also be used to simplify the proof of theorem 4.3.5 in~\citep{book:hiaipetz}.
As in~\citep{book:hiaipetz}, theorem~\ref{432} is sufficient to prove asymptotic freeness almost everywhere for the family
\[
  \left( \left( \{ U(s,n),U(s,n)^{\ast} \} \right)_{s\in S} , \{ D(t,n),D(t,n)^{\ast} : t\in T \} \right)
\]
when the $D_r(n)$ is known to have a limit distribution. It will
also be useful to us that theorem~\ref{432} gives us bounds also
in cases where the $D_r(n)$ do {\em not} converge to a limit. The
$D_r(n)$ model in our case concerns $\frac{1}{\sqrt{n}}R_n$ random
matrices, for which it is not known whether an almost sure limit
exists (only that $\Gamma_n=\frac{1}{N}R_nR_n^{\ast}$ has an
almost sure limit). Also, theorem~\ref{432} gives us grounds for
proving that only the lower mixed moments converge to zero. It can
be applied to cases where only the lower $\| \cdot \|_p$-norms are
known to be bounded, in which only lower mixed moments can be
bounded.

What we really want is to use random matrices $R_n$ independent from the $U_n$ instead of the deterministic matrices $D_r(n)$.
This is addressed by the following theorem.
We restrict to the case of one standard unitary random matrix. 
\begin{teo} \label{23101}
  Let $U_n$ be $n\times n$ standard unitary random matrices, and let $R_n$ be random matrices independent from $U_n$,
  such that $R_nR_n^{\ast}$ converges in distribution almost surely to a compactly supported probability measure $\rho$.
  Then
  \begin{equation} \label{traceconv}
      | tr_n\left( U_n^{m_1}P_1(R_n)U_n^{m_2}P_2(R_n)\cdots U_n^{m_l}P_l(R_n) \right) | \rightarrow 0 \mbox{ a.s.}
  \end{equation}
  uniformly for any choice of polynomials $P_1,...,P_l$ such that $tr_n(P_i(R_n)) = 0$ for all $1\leq i\leq l$.
\end{teo}
The proof is in section~\ref{proof2310}.
As for theorem~\ref{432}, theorem~\ref{23101} is sufficient to prove asymptotic freeness almost everywhere for the family $(U_n,R_n)$ when the $R_n$ are additionally known to have a limit distribution.

The proof is split in two:
First (\ref{traceconv}) is shown for random matrices satisfying bounds of the form $\| R_n \|_p \leq R_p$ ($p\geq 1$).
The proof in this case uses theorem~\ref{432} and is quite short.
The more general case of compactly supported probability measures is proved with an approximation argument.

The next step is to pass from standard unitary random matrices $U_n$ to standard Gaussian random matrices $X_n$.
Note that it could be possible to skip starting with standard unitary random matrices altogether, 
by building directly on results for almost sure convergence of Gaussian random matrices
like those in~\citep{paper:thorbjornsen1}.
We have chosen the approach with standard unitary random matrices for compatibility with~\citep{book:hiaipetz}.
We will prove the following:
\begin{teo} \label{23102}
  Let $X_n$ be $n\times n$ standard Gaussian random matrices, and let $R_n$ be random matrices independent from $X_n$,
  such that $R_nR_n^{\ast}$ converges in distribution almost surely to a compactly supported probability measure $\rho$.
  Then
  \[
    | tr_n\left( Q_1(X_n)P_1(R_n)Q_2(X_n)P_2(R_n)\cdots Q_l(X_n)P_l(R_n) \right) | \rightarrow 0 \mbox{ a.s.}
  \]
  uniformly for any choice of polynomials $Q_q,...,Q_l$, $P_1,...,P_l$ such that 
  \[
    tr_n(Q_i(R_n)) = 0 \mbox{ and } tr_n(P_i(R_n)) = 0
  \]
  for all $1\leq i\leq l$.
\end{teo}
The proof is quite short, and also presented in
section~\ref{proof2310}. Note that the approximation argument used
in the proof of theorem 4.3.5 in~\citep{book:hiaipetz} does not
work in this case. As for theorem~\ref{23101}, theorem~\ref{23102}
is enough to prove asymptotic freeness almost everywhere when the
$R_n$ are additionally known to have a limit distribution. Just as
theorem~\ref{432} gives bounds for mixed moments also in cases
where the deterministic matrices do not converge in distribution,
theorem~\ref{23101} and it's counterpart for Gaussian random
matrices can be used to bound mixed moments in cases where it is
only known that the $R_n$ matrices satisfy $\| \cdot \|_p$-norm
bounds.

To finish the proof we will model our situation through the following theorem, which is stated independently of a random matrix setting.
\begin{teo} \label{prop2}
  Suppose that $a$ and $\{ p,b\}$ are $\ast$-free, with $a$ $R$-diagonal and $p$ a projection with $\phi(p) = c$.
  In the reduced probability space $(p{\cal A}p, \phi(p)^{-1}\phi)$, $\mu_{p(a+b)(a+b)^{\ast}p}$ is uniquely identified by
  $\mu_{paa^{\ast}p}$ and $\mu_{pbb^{\ast}p}$ through the equation
  \begin{equation} \label{connectionrect}
    \mu_{p(a+b)(a+b)^{\ast}p} \boxslash \mu_c = \left( \mu_{pa a^{\ast}p} \boxslash \mu_c \right) \boxplus \left( \mu_{pb b^{\ast}p} \boxslash \mu_c \right)
  \end{equation}
  In particular, $\mu_{p(a+b)(a+b)^{\ast}p}$ has no dependence on mixed moments of $a$ and $b$.
\end{teo}
This will be proved in section~\ref{proofprop2}.
Note that there is no assumption on freeness between $p$ and $b$.
The case $c=1$ is particularly interesting, and corresponds to $p=I$.
In this case, $\mu_{(a+b)(a+b)^{\ast}}$ is uniquely identified by $\mu_{aa^{\ast}}$ and $\mu_{bb^{\ast}}$,
and (\ref{connectionrect}) is simply
\begin{equation} \label{connection}
  \mu_{(a+b)(a+b)^{\ast}} \boxslash \mu_1 = \left( \mu_{a a^{\ast}} \boxslash \mu_1 \right) \boxplus \left( \mu_{b b^{\ast}} \boxslash \mu_1 \right).
\end{equation}
This equation has an interpretation in terms of square random matrices.

Due to (\ref{connectionrect}), $R$-diagonality relieves us from dependencies of many mixed moments, so that some cancellation phenomenon must occur.
This also happens in other cases.
If $a$ and $b$ are free,~\citep{book:comblect} expresses the distribution of the {\em free commutator}, i.e.
\begin{equation} \label{freecommresult}
  R_{\mu_{i(ab-ba)}}(z) = 2\left( R_{\mu_a}^{\mbox{even}} \freestar R_{\mu_b}^{\mbox{even}} \freestar Zeta \right) (z^2),
\end{equation}
where $R^{\mbox{even}}(z) = \sum_{n=1}^{\infty} \alpha_{2n} z^n$ whenever $R(z) = \sum_{n=1}^{\infty} \alpha_n z^n$.
(\ref{freecommresult}) holds also when $a$ and $b$ are not $R$-diagonal.
(\ref{freecommresult}) also expresses a connection with multiplicative free convolution with $\mu_1$,
since boxed convolution with the $Zeta$-series is involved.

Theorem~\ref{prop2} has a more general flavour than theorem~\ref{teo1a},
since the limits $\frac{1}{N} X_nX_n^{\ast}$ from theorem~\ref{teo1a} do not include all $R$-diagonal pairs.
The following limiting version of theorem~\ref{prop2} will be useful in finishing the proof of theorem~\ref{teo1a}:

\begin{teo} \label{prop2b}
  Let the random variables $\{ a_n,b_n,p_n \}\in ({\cal A}_n,\phi_n)$, $\{ a,p \}\in ({\cal A},\phi)$ be given,
  where $a$ is $R$-diagonal and $p,p_n$ are projections with $\phi(p) = \phi_n(p_n) = c$.
  Form the random variable $paa^{\ast}p$ in $( p{\cal A}p , \phi(p)^{-1}\phi)$,
  and the random variables $p_nb_nb_n^{\ast}p_n$ and $p_n(a_n+b_n)(a_n+b_n)^{\ast}p_n$ in $(p_n{\cal A}_n p_n,\phi(p_n)^{-1}\phi_n)$.
  If
  \[
    \mu_{a_n,a_n^{\ast}} \rightarrow \mu_{a,a^{\ast}} \mbox{, } \mu_{p_n} \rightarrow \mu_p,
  \]
  and
  \[
    \mu_{p_nb_nb_n^{\ast}p_n} \rightarrow \mu,
  \]
  in distribution,
  moments are uniformly bounded in $n$, and mixed moments of $(a_n,\{ p_n,b_n \})$ go to $0$,
  then $\mu_{p_n(a_n+b_n)(a_n+b_n)^{\ast}p_n}$ converges in distribution and the limit is uniquely identified by the equation
  \begin{equation} \label{inthelimit}
    \lim_{n\rightarrow\infty} \mu_{p_n(a_n+b_n)(a_n+b_n)^{\ast}p_n} \boxslash \mu_c =
    \left( \mu_{paa^{\ast}p} \boxslash \mu_c \right) \boxplus \left( \mu \boxslash \mu_c \right)
  \end{equation}
\end{teo}
\begin{proof}
The limiting moments of $\mu_{p_n(a_n+b_n)(a_n+b_n)^{\ast}p_n} \boxslash \mu_c$
do not change if we "zero out" the mentioned mixed moments (i.e. that we assume freeness of $( a_n, \{ p_n,b_n\})$),
due to the assumption on their vanishing and of uniform boundedness on moments.
It is also easily seen that the limiting moments do not change if we change the distribution of $a_n$ to $\mu_{a_n} = \mu_a$ for all $n$. But then
\[
  \lim_{n\rightarrow\infty} \mu_{p_n(a_n+b_n)(a_n+b_n)^{\ast}p_n} \boxslash \mu_c =
  \lim_{n\rightarrow\infty} \mu_{p_n(a+b_n)(a+b_n)^{\ast}p_n} \boxslash \mu_c
\]
\[
  = \lim_{n\rightarrow\infty} \left( \mu_{p_na a^{\ast}p_n} \boxslash \mu_c \right) \boxplus \left( \mu_{p_nb_n b_n^{\ast}p_n} \boxslash \mu_c \right) =
  \left( \mu_{paa^{\ast}p} \boxslash \mu_c \right) \boxplus \left( \mu \boxslash \mu_c \right)
\]
where we used theorem~\ref{prop2}, so that (\ref{inthelimit}) holds.
\end{proof}

The rest of the proof of theorem~\ref{teo1a} now goes as follows:
The rectangular random matrices $R_n$ can be viewed as the $N\times N$ random matrices $p_nS_n$,
where the projection $p_n$ is a diagonal constant matrix, with the fraction of $1$'s on the diagonal equal to $c$,
and $S_n$ is an extension of the $n\times N$ matrix $R_n$ to an $N\times N$-matrix, obtained by adding zeros.
Similarly, the random matrices $X_n$ can be viewed as the $N\times N$-matrices $p_nY_n$,
where $Y_n$ is an extension of the $n\times N$ matrix $X_n$ to an $N\times N$-matrix, obtained by adding more independent standard Gaussian entries.

Since $\frac{1}{N}S_nS_n^{\ast}$ almost surely converges to a compactly supported probability measure,
$(\frac{1}{\sqrt{n}} Y_n , \frac{1}{\sqrt{n}} S_n)$ satisfies the requirements of theorem~\ref{23102}.
Thus, mixed moments of $\frac{1}{\sqrt{n}} X_n$ and $\frac{1}{\sqrt{n}} S_n$ go to zero almost surely.
It is also seen that $\frac{1}{\sqrt{n}} S_n$ has it's moments bounded as $n\rightarrow\infty$ almost surely.
It is well known~\citep{book:hiaipetz} that $\frac{1}{\sqrt{N}}Y_n$ converges in distribution almost surely to the circular law, which is $R$-diagonal.

Thus, all assumptions of theorem~\ref{prop2b} are satisfied for $a_n = \frac{1}{\sqrt{N}} Y_n$, $b_n = \frac{1}{\sqrt{N}} S_n$ and $p_n$, almost surely.
Thus, almost surely,
\[
  \lim_{n\rightarrow\infty} \mu_{p_n\frac{1}{N} (S_n+\sigma Y_n)(S_n+\sigma Y_n)^{\ast}p_n} \boxslash \mu_c =
  \left( \mu_{p \sigma^2 \frac{1}{N}YY^{\ast}p} \boxslash \mu_c \right) \boxplus \left( \mu_{\Gamma} \boxslash \mu_c \right)
\]
\[
  = \left( \mu_{\sigma^2 I} \boxtimes \mu_c \boxslash \mu_c \right) \boxplus \left( \mu_{\Gamma} \boxslash \mu_c \right) =
  \mu_{\sigma^2 I} \boxplus \left( \mu_{\Gamma} \boxslash \mu_c \right),
\]
or
\[
  \lim_{n\rightarrow\infty} \mu_{\frac{1}{N} (R_n+\sigma X_n)(R_n+\sigma X_n)^{\ast}} \boxslash \mu_c = \mu_{\sigma^2 I} \boxplus \left( \mu_{\Gamma} \boxslash \mu_c \right),
\]
which is the statement of theorem~\ref{teo1a}.

The proof as skecthed here assumes $c\leq 1$. An explanation for how the proof goes for $c > 1$ is given succeeding the proof of theorem~\ref{prop2}.

\subsection{The proof of theorem~\ref{432}} \label{proof432}
The proof will use the (generalized) H\"{o}lder inequality:
\begin{lem} \label{holder}
For matrices $A_1,...,A_k$, the following holds:
\[
  \| A_1 \cdots A_k \|_p \leq \| A_1 \|_{p_1} \cdots \| A_k \|_{p_k} \mbox{ when } \sum_{i=1}^k \frac{1}{p_i} = \frac{1}{p}.
\]
\end{lem}

In the proof of lemma 4.3.2 in~\citep{book:hiaipetz}, (\ref{432expression}) is written as
\[
  \left( \frac{1}{n} \right)^2 \sum_{i_1,...,i_{2k} = 1}^n \sum_{j_{k(1)},...,j_{k(l)},j_{k(l+1)+1},...,j_{k(2l)+1}}^n \left( \prod_{r=1}^l d_{ j_{k(r)} i_{k(r) + 1} }(t_r,n) \right)
\]
\begin{equation} \label{writtenout}
  \times \left( \prod_{r=l+1}^{2l} \bar{d}_{ i_{k(r)} j_{k(r) + 1} }(t_r,n) \right) E \left( \prod_{h=1}^{2k} u_{i_h j_h} (s(h),\varepsilon(h),n) \right),
\end{equation}
where for $1\leq r\leq l$,
\[
  k(r) = |m_1| + \cdots + |m_r|,
\]
\[
  k=k(l) \mbox{, } k(l+r) = k + k(r) \mbox{, } t_{l+r} = t_r.
\]
Moreover, for $h$ such that $k(r-1) + 1 \leq h \leq k(r)$,
\[
  s(h) = s_r \mbox{, } \varepsilon(h) =
  \left\{
    \begin{array}{ll}
      1 & \mbox{ if } m_r>0 \\
      -1 & \mbox{ if } m_r<0
    \end{array}
  \right.
\]
Here $s(h+k) = s(h)$ and $\varepsilon(k+h) = -\varepsilon(h)$ for $1\leq h\leq k$, and
\[
  u_{ij}(s,\varepsilon ,n) = \left\{ \begin{array}{ll} U_{ij}(s,n) & \mbox{ if } \varepsilon = 1 \\ \bar{U}_{ji}(s,n) & \mbox{ if } \varepsilon = -1 \end{array} \right.
\]
Since (\ref{writtenout}) is a matrix product written out, the following must hold:
\begin{equation} \label{ekvrel1}
  \left\{ \begin{array}{ll} j_h = i_{h+1} & \mbox{ for } h\in \{ 1,...,k \} \setminus \{ k(1),...,k(l)\} ,\\ i_h = j_{h+1} & \mbox{ for } h\in \{ k+1,...,2k\}\setminus\{ k(l+1),...,k(2l)\} \end{array} \right.
\end{equation}
Also, due to the vanishing of many mixed moments of entries in standard unitary random matrices (lemma 4.2.2 in~\citep{book:hiaipetz}),
two pair partitions ${\cal U}$ and ${\cal V}$ can be chosen so that if $\{ h,h' \}\in {\cal U}$ then
\begin{equation} \label{ekvrel2}
  s(h) = s(h') \mbox{, } \varepsilon(h) = 1 \mbox{, } \varepsilon(h') = -1 \mbox{, } i_h = j_{h'} (=i_{h'}+1),
\end{equation}
and if $\{ h,h' \}\in {\cal V}$ then
\begin{equation} \label{ekvrel3}
  s(h) = s(h') \mbox{, } \varepsilon(h) = -1 \mbox{, } \varepsilon(h') = 1 \mbox{, } i_h = j_{h'} (=i_{h'}+1)
\end{equation}
These two pair partitions and (\ref{ekvrel1}) cause many equalities among the $i_1,...,i_{2k}$,
and define the equivalence relation ${\cal R(U,V)}$ on $\{ 1,...,2k\}$ so that $i_h = i_{h'}$ whenever $h$ and $ h'$ are in the same equivalence class of ${\cal R(U,V)}$.
We let $k_0$ denote the number of equivalence classes of ${\cal R(U,V)}$, and let $h(1),...,h(k_0)$ be representatives from the equivalence classes.

Recall the expressions
\begin{equation} \label{ciota}
  C_n(\iota_1,...,\iota_{k_0}) = \left( \prod_{r=1}^l d_{ j_{k(r)} i_{k(r)+1} }(t_r,n) \right) \left( \prod_{r=l+1}^{2l} \bar{d}_{ i_{k(r)} j_{k(r)+1} }(t_r,n)  \right),
\end{equation}
\begin{equation} \label{qnexpression}
  Q_n(\iota_1,...,\iota_{k_0}) = E\left( \prod_{h=1}^{2k} u_{i_h j_h}(s(h),\varepsilon(h),n) \right)
\end{equation}
from~\citep{book:hiaipetz},
where $(\iota_1,...,\iota_{k_0})$ in (\ref{ciota}) are defined as $(i_{h(1)},...,i_{h(k_0)})$ (i.e. representatives of the equivalence classes),
and $j_1,...,j_{2k}$ are determined subject to ${\cal U,V}$.
(4.3.6) of~\citep{book:hiaipetz} says that it is enough to prove that for any partition ${\cal W}$ of $\{ 1,...,k_0\}$ we have that
\begin{equation} \label{436bound}
  \sum_{(\iota_1,...,\iota_{k_0}):{\cal W}} C_n(\iota_1,...,\iota_{k_0}) Q_n(\iota_1,...,\iota_{k_0}) = O(1) \mbox{ as } n\rightarrow\infty ,
\end{equation}
where the summation is over $(\iota_1,...,\iota_{k_0})$ such that $\iota_p = \iota_q$ if and only if $p$ and $q$ are in the same block of ${\cal W}$.
For a given ${\cal W}$,
it is known that $| Q_n(\iota_1,...,\iota_{k_0}) | = O(n^{-k})$ uniformly for $i_k,j_k$ as $n\rightarrow\infty$,
and that it has the same value for all $(\iota_1,...,\iota_{k_0})$ taking part in the sum (\ref{436bound}).
So, from (\ref{startsum}) we deduce that it is enough to show that
for any choice of ${\cal U,V,W}$ (there is a finite number of such choices),
\begin{equation} \label{startsum}
  \sum_{(\iota_1,...,\iota_{k_0}):{\cal W}} C_n(\iota_1,...,\iota_{k_0})
\end{equation}
is $O(n^k)$. In\citep{book:hiaipetz} this is proved using the fact that (\ref{ciota}) are bounded uniformly.
This is not true in our case since only uniform boundedness in $\| \cdot \|_p$ is assumed.
Instead, we will group sums of terms into matrix multiplication units, and use the H\"{o}lder inequality
together with the $\| \cdot \|_p$-norm bounds.
Instead of the terms in (\ref{startsum}), where the sum is over
\[
  (\iota_1,...,\iota_{k_0}) \mbox{: } \iota_p = \iota_q \mbox{ if and only if $p$ and $q$ are in the same block of ${\cal W}$,}
\]
it will be better for us to sum over
\[
  (\iota_1,...,\iota_{k_0}) \mbox{: }  \iota_p = \iota_q \mbox{ if $p$ and $q$ are in the same block of ${\cal W}$.}
\]
The latter set is more compatible with indices in multiplications of many matrices.
This second set is larger than the first, and can be written as
\begin{equation} \label{newsum}
  \sum_{{\cal W'\geq W}} \sum_{(\iota_1,...,\iota_{k_0}):{\cal W'}} C_n(\iota_1,...,\iota_{k_0})
\end{equation}
It is obvious that (\ref{startsum}) can be written
\begin{equation} \label{splitsum}
   \sum_{{\cal W'\geq W}} a_{{\cal W'}} \sum_{{\cal W''\geq W'}} \sum_{(\iota_1,...,\iota_{k_0}):{\cal W''}} C_n(\iota_1,...,\iota_{k_0}),
\end{equation}
where $a_{{\cal W'}}$ are integer constants which can easily be calculated
(proving (\ref{splitsum}) boils down to splitting all values of $\iota_1,\iota_2$ into those where $\iota_1 = \iota_2$, and those where $\iota_1\neq\iota_2$.
This is done recursively and for all $\iota_i$ to yield (\ref{splitsum})).
Since there is a finite number of elements in the two outer sums in (\ref{splitsum}), to prove that (\ref{startsum}) is $O(n^k)$
it is enough to show that (\ref{newsum}) is $O(n^k)$ for any choice of ${\cal W}$.

Let $l_0$ denote the number of equivalence classes in ${\cal R(U,V)}$ with only one entry,
and let $h(1),...,h(l_0)$ be the corresponding respresentatives.
According to~\citep{book:hiaipetz}, equivalence classes with only one entry give rise to factors of the form $\tilde{d}_{\iota_i \iota_i}$ in (\ref{ciota}),
where $\tilde{d}$ is either $d_{\iota_i \iota_i}(t_r,n)$ or $\bar{d}_{\iota_i \iota_i}(t_r,n)$ for some $r$.
Equivalence classes with only one entry thus leads (through summation over one $\iota_i$ appearing in just one factor)
to factors in (\ref{newsum}) which are (non-normalized) traces of the $D(t_r,n)$.
These are zero, so we can assume that $l_0 = 0$ when we attempt to bound (\ref{newsum}).
Had we used the sum (\ref{startsum}) instead of (\ref{newsum}), we would not obtain zero.

So we assume that there are no singleton equivalence classes, i.e. $k_0\leq k$.
Let $K_0$ be the number of equivalence classes actually appearing in (\ref{ciota}) (this is a function of ${\cal U}$ and ${\cal V}$).
we have that $K_0 \leq k_0$, but equality does not necessarily hold.
We will use matrix units $E_{ij}$ (i.e. $E_{ij}(i,j) = \delta_{ij}$).
By placing matrix units $F_i$ with indices from $\iota_1,...,\iota_{k_0}$ in between the terms in (\ref{ciota}),
(\ref{newsum}) can be written as
\begin{equation} \label{torearrange}
  \sum_{\cal W'\geq W} \sum_{(\iota_1,...,\iota_{k_0}):{\cal W'}} ntr_n\left( \prod_{i=1}^{2l}(F_i D_i)\right),
\end{equation}
where $D_i$ are matrices from $D(t_r,n)$ or one of their transposes/conjugates.
Since $\| E_{ij} \|_p = n^{-\frac{1}{p}}$ and the number of possible choices of matrix units is $n^{K_0}$, lemma~\ref{holder} implies that
(\ref{torearrange}) is bounded by
\begin{equation} \label{torearrange2}
  n^{K_0+1} n^{-\frac{1}{2}} \prod_{i=1}^{2l} \| D_i \|_{4l} = n^{K_0+\frac{1}{2}} \prod_{i=1}^{2l} \| D_i \|_{4l}.
\end{equation}
Since $K_0\leq k$, this is $O(n^k)$ except possibly in the case when $K_0=k$,
i.e. when all equivalence classes have exactly two elements.

So, for the rest of the proof, we assume that all equivalence classes have exactly two elements.
Note that the number of times an equivalence class appears as an $i$ is equal to the number of the times the same class appears as a $j$ in (\ref{ciota}).
This is obvious from the way the equivalence relation is defined (\ref{ekvrel1}), (\ref{ekvrel2}), (\ref{ekvrel3})
in order to avoid a zero value in (\ref{qnexpression}).
This means that we can take the first of the $K_0$ equivalence classes appearing in (\ref{ciota}), 
and rearrange the terms in (\ref{ciota}) so that
the equivalence class appear in alternating order as an $i$ and as a $j$. (\ref{torearrange}) can thus be rewritten to
\begin{equation} \label{torearrange3}
  \sum_{\cal W'\geq W} \sum_{(\iota_1,...,\iota_{k_0}):{\cal W'}} ntr_n \left( F_1 D_1 G_1 D_2 F_2 D_3 G_2 D_4 F_3 \right),
\end{equation}
where $D_1,D_2,D_3,D_4$ are the matrices where the first equivalence class appear as an $i$ or a $j$, and in alternating order. 
Also, $F_1=F_2=F_3=E_{rr}$ are matrix units, $r$ is a given number, and the $G_i$ are products man of the matrices $D_i$ in (\ref{torearrange}).
(\ref{torearrange3}) can also be written
\begin{equation} \label{torearrange4}
  \sum_{\cal W'\geq W} \sum_{(\iota_1,...,\iota_{k_0}):{\cal W'}} ntr_n \left( diag(D_1 G_1 D_2) diag(D_3 G_2 D_4) \right),
\end{equation}
where $diag(A)$ stands for the diagonal of the matrix $A$.
Similarly to the calculation of the bound (\ref{torearrange2}), (\ref{torearrange4}) is seen to be bounded by
\begin{equation} \label{torearrange5}
  n^{K_0} \| diag(D_1 G_1 D_2) \|_2 \| diag(D_3 G_2 D_4) \|_2.
\end{equation}
Note that $\| diag(A) \|_2 \leq \| A \|_2$, since
\[
  \| diag(A) \|_2^2 = \frac{1}{n} \sum_i |a_{ii}|^2 \leq \frac{1}{n} \sum_{i,j} |a_{ij}|^2 = tr_n(A^{\ast}A) = \| A \|_2^2,
\]
where $A=(a_{ij})_{i,j}$.
This means that (\ref{torearrange5}) is bounded by
\[
  n^{K_0} \| D_1 G_1 D_2 \|_2 \| D_3 G_2 D_4 \|_2,
\]
which is $O(n^{K_0})$ and hence $O(n^k)$ since all $D_i$ are bounded in $p$-norm, and the only other factors are matrix units, which
have $p$-norm $n^{-\frac{1}{p}}$. 

That there exists a $p_l$ for a given $l$ as in the last statement of the theorem is obvious from the proof and the 
way the H\"{o}lder inequality was used. This completes the proof.

\subsection{The proofs of theorem~\ref{23101} and~\ref{23102}} \label{proof2310}
First assume that $R_n$ satisfies $\| R_n \|_p \leq R_p$ ($p\geq 1$) almost surely for some constants $R_p$.
$P_i(R_n)$ satisfies similar $\| \cdot \|_p$-norm bounds due to lemma~\ref{holder}.
Call the underlying probability space $\Omega$.
Denote by $f_{U_n,R_n}(U,R)$ the joint density of $R_n$ and $U_n$, and by $f_{U_n}(U)$ and $f_{R_n}(R)$ the marginal densities.
Due to independence, $f_{U_n,R_n}(U,R) = f_{U_n}(U) f_{R_n}(R)$, and therefore
\[
  E\left( | tr_n\left( U_n^{m_1} P_1(R_n) U_n^{m_2} P_2(R_n) \cdots U_n^{m_l} P_l(R_n) \right) |^2 \right)
\]
\begin{equation} \label{independenceequation}
  \begin{array}{ll}
    =& \int_{\Omega} | tr_n\left( U_n^{m_1} P_1(R_n) \cdots U_n^{m_l} P_l(R_n) \right) |^2 ds\\
    =& \int_{M_n(\C)} \int_{M_n(\C)}   | tr_n\left( U^{m_1} P_1(R) \cdots U^{m_l} P_l(R) \right) |^2 f_{U_n,R_n}(U,R) dU dR\\
    =& \int_{M_n(\C)} \int_{M_n(\C)}   | tr_n\left( U^{m_1} P_1(R) \cdots U^{m_l} P_l(R) \right) |^2 f_{U_n}(U) dU f_{R_n}(R) dR \\
\leq & \int_{M_n(\C)} Cn^{-2} f_{R_n}(R) dR = Cn^{-2}
  \end{array}
\end{equation}
where we have used the bounds for deterministic matrices from theorem~\ref{432}.
Therefore
\[
  | tr_n\left( U_n^{m_1}P_1(R_n)U_n^{m_2}P_2(R_n)\cdots U_n^{m_l}P_l(R_n) \right) | \rightarrow 0 \mbox{ a.s.}
\]
for such random matrices $R_n$. If $R_nR_n^{\ast}$ is just known to converge in distribution almost surely to a compactly supported probability measure,
observe that almost surely there exists a value $R$ so that $\| R_n \|_p \leq R$ for $n$ large enough~\citep{book:hiaipetz}.
For each $l$, choose $p_l$ as in the statement of theorem~\ref{432}.
Denote by $\Omega_{p_l,N}$, $p\geq 1$ $N\in\N$ the subset of $\Omega$ determined by values $s$ such that
\begin{equation}
  \| R_n(s) \|_{p_l} \leq N
\end{equation}
for all $n$. Define $R_{n,p_l,N} = \chi_{\Omega_{p_l,N}}R_n$ with
$\chi$ denoting the characteristic function. The $R_{n,p_l,N}$
satisfy the estimates (\ref{independenceequation}) for mixed
moments of length $\leq l$, so that these mixed moments go to zero
almost surely in $\Omega_{p_l,N}$. $\cup_{p_l,N} \Omega_{p_l,N}$
has probability 1: Almost surely, the $\| \cdot \|_{p_l}$-norm of
$R_n$ stays bounded by some finite value for large enough $n$.
Thus, for every $s$ in a set with probability one, we can find a
value $N_s$ such that $\| R_n(s) \|_{p_l} \leq N_s$ for ALL $n$.
But then $s\in \Omega_{p_l,N_s}$, so that $\cup_{p_l,N}
\Omega_{p_l,N}$ has probability 1 as claimed. Since
$\cup_N\Omega_{p_l,N}$ has probability 1, theorem~\ref{23101}
follows from the fact that $R_n = R_{n,p_l,N}$ on
$\Omega_{p_l,N}$. By increasing $l$ we get almost sure convergence
of higher mixed moments to zero also.

Now for theorem~\ref{23102}. Write
\[
  X_n = U_n \Lambda_n U_n^{\ast}
\]
for a unitary random matrix $U_n$, and diagonal random matrix $\Lambda_n$. 
We may assume that $U_n$ is a
standard unitary random matrix, as in the proof of theorem 4.3.5
of~\citep{book:hiaipetz}, since Gaussian random matrices are
unitarily invariant. We can also assume that $\Lambda_n$ is
independent from $U_n$, so that $(U_n,\{ \Lambda_n , R_n
\})$ is an independent family. 
$R_n$ converges to a limit which is compactly supported, 
and $\Lambda_n$ does the same. 
Since (\ref{independenceequation}) can be easily generalized to the case where the $R_n$ 
are replaced with many different $R_n^î$ (with the $\{ R_n^î \}$ all independent from $U_n$) 
we conclude also for theorem~\ref{23102} that we get
almost sure convergence to zero of mixed moments as in definition~\ref{freedef}.

\subsection{The proof of theorem~\ref{prop2}} \label{proofprop2}
First write $\phi\left( (p(a+b)(a+b)^{\ast}p)^m \right)$ as a sum of mixed moments of length $3m$ by multiplying out $(p(a+b)(a+b)^{\ast}p)^m$:
\begin{equation} \label{multipliedout}
  \phi\left( (p(a+b)(a+b)^{\ast}p)^m \right) = \sum_{ \sigma_1\leq\sigma } \phi\left( x_1 x_2^{\ast} p \cdots x_{2m-1} x_{2m}^{\ast} p \right),
\end{equation}
where $\sigma = \{ 1,2, 4,5,\cdots 3m-2,3m-1 \}$
($\{ 1,2, 4,5,\cdots 3m-2,3m-1 \}$ correspond to the indices of the locations of the $x_i$, $x_i^{\ast}$
in the moments $x_1 x_2^{\ast} p \cdots x_{2m-1} x_{2m}^{\ast} p$),
$\sigma_1$ runs over all subsets of $\sigma$,
and $x_i = a$ if $i\in\sigma_1$, $x_i = b$ if $i\in\sigma\setminus\sigma_1$.
Denote by $|\sigma_1|$ the cardinality of $\sigma_1$.
We denote by $\alpha$ the cumulants of $\mu_{a,a^{\ast}}$ and $\beta$ the cumulants of $\mu_{b,b^{\ast},p}$,
so that the moment-cumulant formula for $(a,a^{\ast})$ is
\begin{equation} \label{cum1}
  \phi\left( x_{i_1} \cdots x_{i_n} \right) = \sum_{\pi = \{ B_1 , \cdots , B_k \} \in NC(n)} \prod_{i=1}^k [ coef(i_1,...,i_n)|B_i ](R_{\mu_{a,a^{\ast}}})
\end{equation}
with $x_1=a$, $x_2=a^{\ast}$ and $i=1$ or $2$,
and the moment-cumulant formula for $(b,b^{\ast},p)$ is
\begin{equation} \label{cum2}
  \phi\left( x_{i_1} \cdots x_{i_n} \right) = \sum_{\pi = \{ B_1 , \cdots , B_k \} \in NC(n)} \prod_{i=1}^k [ coef(i_1,...,i_n)|B_i ](R_{\mu_{b,b^{\ast},p}}),
\end{equation}
with $x_1=b$, $x_2=b^{\ast}$, $x_3=p$ and $i=1$, $2$ or $3$.
We will use the shorthand notation
\[
  \begin{array}{lll}
    \alpha_{B_i} &=& [ coef(i_1,...,i_n)|B_i ](R_{\mu_{a,a^{\ast}}})\\
    \beta_{B_i}  &=& [ coef(i_1,...,i_n)|B_i ](R_{\mu_{b,b^{\ast},p}})
  \end{array}
\]
Due to the freeness of $a$ and $\{ p,b\}$, the moment-cumulant formula applied to all moments in (\ref{multipliedout}) and (\ref{usefulsum}) yields
\begin{equation} \label{sumnoncrossing}
  \sum_{ \sigma_1\leq\sigma }
  \sum_{ \stackrel{\pi_1\in NC(|\sigma_1|)}{\pi_1\leq\sigma_1} }
  \sum_{ \stackrel{\pi_2\in NC(|\sigma_1^c|)}{\pi_2\leq\sigma_1^c \mbox{, no crossings between $\pi_1$ and $\pi_2$}} }
  \alpha_{\pi_1} \beta_{\pi_2},
\end{equation}
where $\sigma_1^c= \{ 1,...,3m\}\setminus\sigma_1$ and $|\sigma_1^c|$ is the cardinality of $\sigma_1^c$.
$\pi_1$ divides $\{ 1,...,|\sigma_1|\}$ into $|K(\pi_1)|$ sets
(see definition~\ref{krewerasdef} and figure~\ref{fig:big})
when $\pi_1$ is viewed as an element in $NC(|\sigma_1|)$.
$\pi_1$ also divides $\{ 1,...,3m\}$ into the same number of sets,
according to the circular representation of $\{ 1,...,3m\}$.
Let us denote these blocks by $B_1,\cdots B_k$, so that
$\sigma_1^c = \{ B_1 , \cdots , B_k \}$ as a subpartition of $1_{|\sigma_1^c|}$.
Since $\pi_1$ and $\pi_2$ have no crossings if and only if $\pi_2\leq \{ B_1,...,B_k\}$,
(\ref{sumnoncrossing}) can be written as
\begin{equation} \label{sumnoncrossing1}
  \sum_{ \sigma_1\leq\sigma }
  \sum_{ \stackrel{\pi_1\in NC(|\sigma_1|)}{\pi_1\leq\sigma_1} }
  \sum_{ \stackrel{\pi_2\in NC(|\sigma_1^c|)}{\pi_2\leq \{ B_1 , \cdots , B_k \} } }
  \alpha_{\pi_1} \beta_{\pi_2}
\end{equation}
When $\pi_2\leq \{ B_1,...,B_k\}$ we can write $B_i = \pi_{2i1} \cup \pi_{2i2}, \cdots$,
where the $\pi_{2ij}$ are the reindexed blocks of $\pi_2$ which are contained in $B_i$,
and where $\pi_{2i} = \{\pi_{2i1} , \pi_{2i2}, \cdots \}$.
This is in $NC(|B_i|)$ since $\pi_2$ is noncrossing.
First rewrite (\ref{sumnoncrossing1}) to
\begin{equation}
  \sum_{ \sigma_1\leq\sigma }
  \sum_{ \stackrel{\pi_1\in NC(|\sigma_1|)}{\pi_1\leq\sigma_1} }
  \alpha_{\pi_1}
  \left( \sum_{ \pi_{2i}\in NC(|B_i|)} \prod_{i=1}^k \beta_{\pi_{2i}} \right).
\end{equation}
Then note that the $\pi_{2i}\in NC(|B_i|)$ can be summed independently of one another, so that we can rewrite to
\begin{equation} \label{sumnoncrossing2}
  \sum_{ \sigma_1\leq\sigma }
  \sum_{ \stackrel{\pi_1\in NC(|\sigma_1|)}{\pi_1\leq\sigma_1} }
  \alpha_{\pi_1}
  \prod_{i=1}^k \left( \sum_{ \pi_{2i}\in NC(|B_i|)} \beta_{\pi_{2i}} \right).
\end{equation}
Note also that only $\pi_1$ with blocks 
\[
  C_k = \{ c_{k1},...,c_{kr} \}
\]
where $x_{i_{c_{k1}}},...,x_{i_{c_{kr}}}$ are alternating values of $a$ and $a^{\ast}$,  
give contribution in (\ref{sumnoncrossing2}), due to $R$-diagonality of $a$. 
Hold such a $\pi_1$ fixed in (\ref{sumnoncrossing2}), 
and take a look at the inner sum in (\ref{sumnoncrossing2}) for a given $i$. 
This is simply the moment-cumulant formula (\ref{cum2}) for a moment of length $|B_i|$, where the mixed moment is on the form 
\[
  p             b b^{\ast} p    b b^{\ast} p    \cdots    b b^{\ast} p,
\]
or on the form
\[
  b^{\ast} p    b b^{\ast} p    b b^{\ast} p    \cdots    b b^{\ast} p      b
\]
due to the alternating structure in (\ref{multipliedout}). 
In both cases the moment-cumulant formula yields $\phi\left( (pbb^{\ast}p)^{\frac{|B_i|}{2}} \right)$ 
for the inner sum in (\ref{sumnoncrossing2}). 
Therefore, we get that (\ref{sumnoncrossing2}) equals
\begin{equation} \label{sumnoncrossing3}
  \sum_{ \sigma_1\leq\sigma }
  \sum_{ \stackrel{\pi_1\in NC(|\sigma_1|)}{\pi_1\leq\sigma_1} }
  \alpha_{\pi_1}
  \prod_{i=1}^k
  \phi \left( (pbb^{\ast}p)^{\frac{|B_i|}{2}} \right).
\end{equation}
Since $a$ is $R$-diagonal, the $\alpha_{\pi_1}$ which give contribution in (\ref{sumnoncrossing3}) are uniquely identified by the moments $\phi((aa^{\ast})^m)$.
Therefore, the moments 
\[
  \phi\left( (p(a+b)(a+b)^{\ast}p)^m \right)
\]
are entirely identified by the moments $\phi((aa^{\ast})^m)$ and $\phi((pbb^{\ast}p)^m)$,
so that $\mu_{p(a+b)(a+b)^{\ast}p}$ only depends on $\mu_{a a^{\ast}}$ and $\mu_{pb b^{\ast}p}$.
All distributions are here in the noncommutative probability space $({\cal A},\phi)$, not yet in the reduced space $(p{\cal A}p,\phi(p)^{-1}\phi)$.
If we can prove the theorem when $p$ and $b$ are free, it will also hold when $p$ and $b$ are not free since
$\mu_{p(a+b)(a+b)^{\ast}p}$ only depends on $\mu_{a a^{\ast}}$ and $\mu_{pb b^{\ast}p}$.

We can replace $\mu_{b,b^{\ast}}$ with the (unique) $R$-diagonal pair $b_0$ so that $ \mu_{b b^{\ast}} = \mu_{b_0 b_0^{\ast}}$.
So, we assume that both $a$ and $b$ give rise to free $R$-diagonal pairs.
Their determining series are $R_{\mu_{aa^{\ast}}} \freestar Moeb$ and $R_{\mu_{bb^{\ast}}} \freestar Moeb$, respectively.
$(a+b,(a+b)^{\ast})$ is also an $R$-diagonal pair, with determining series $R_{\mu_{aa^{\ast}}} \freestar Moeb + R_{\mu_{bb^{\ast}}} \freestar Moeb$.
This means that
\begin{equation} \label{connection2}
  R_{\mu_{(a+b)(a+b)^{\ast}}} \freestar  Moeb = R_{\mu_{aa^{\ast}}} \freestar Moeb + R_{\mu_{bb^{\ast}}} \freestar Moeb.
\end{equation}
If $x$ is free from $p$, we next calculate $R_{\mu_{pxx^{\ast}p}}$ in the reduced space $(p{\cal A}p,\phi(p)^{-1}\phi)$.
We call this $R_{\mu_{pxx^{\ast}p}}^{p{\cal A}p}$ in the rest of the proof, with similar notation for the moment series
(In the rest of the paper, this notation is dropped since $pxx^{\ast}p$ is assumed to be in $(p{\cal A}p,\phi(p)^{-1}\phi)$).
Note that $M_{\mu_{pxx^{\ast}p}}^{p{\cal A}p} = \frac{1}{c} M_{\mu_{pxx^{\ast}p}}$.
We have
\[
\begin{array}{lll}
  R_{\mu_{pxx^{\ast}p}}^{p{\cal A}p} &=& M_{\mu_{pxx^{\ast}p}}^{p{\cal A}p} \freestar Moeb = \left( \frac{1}{c}M_{\mu_{pxx^{\ast}p}} \right) \freestar Moeb\\
                                  &=& \left( \frac{1}{c}(R_{\mu_{xx^{\ast}}} \freestar M_{\mu_p}) \right) \freestar Moeb = \left( \frac{1}{c}\left( R_{\mu_{xx^{\ast}}} \freestar (cZeta)\right) \right) \freestar Moeb\\
                                  &=& \left( \frac{1}{c}\left( R_{\mu_{xx^{\ast}}} \freestar (cZeta)\right) \right) \freestar \left( \frac{1}{c}(cMoeb) \right) \\
                                  &=& c^{-n-1}\left( R_{\mu_{xx^{\ast}}} \freestar (cZeta) \freestar (cMoeb) \right)\\
                                  &=& c^{-n-1} \left( R_{\mu_{xx^{\ast}}} \freestar (c^{n+1}Id) \right) = c^{-n-1} \left( R_{\mu_{xx^{\ast}}} \freestar (c^2 Id) \right) \\
                                  &=& c^{-n-1}\left( c^{2n} R_{\mu_{xx^{\ast}}} \right) = c^{n-1} R_{\mu_{xx^{\ast}}}.
\end{array}
\]
For a general Mar\u{c}henko Pastur law $\mu_d$, $R_{\mu_d} = d^{n-1}Zeta$, and it is easily verified that $R_{\mu_d}^{-1} = d^{n-1}Moeb$.
We have that
\[
  R_{\mu_{pxx^{\ast}p}}^{p{\cal A}p} \freestar R_{\mu_d}^{-1} = \left( c^{n-1} R_{\mu_{xx^{\ast}}} \right) \freestar (d^{n-1}Moeb ) = c^n d^n \left( \left( c^{-1} R_{\mu_{xx^{\ast}}} \right) \freestar (d^{-1}Moeb) \right)
\]
if $c=d$, this can be simplified to
\begin{equation} \label{cequalsd}
  c^n c^n c^{-n-1} \left( R_{\mu_{xx^{\ast}}} \freestar Moeb \right) = c^{n-1} \left( R_{\mu_{xx^{\ast}}} \freestar Moeb \right).
\end{equation}
Using this for $x=a$, $x=b$ and $x=a+b$, and also using (\ref{connection2}), we get
\[
\begin{array}{lll}
  R_{\mu_{p(a+b)(a+b)^{\ast}p}}^{p{\cal A}p} \freestar R_{\mu_c}^{-1} &=& c^{n-1} \left( R_{\mu_{(a+b)(a+b)^{\ast}}} \freestar Moeb \right)\\
                                                                   &=& c^{n-1} \left( R_{\mu_{aa^{\ast}}} \freestar Moeb + R_{\mu_{bb^{\ast}}} \freestar Moeb \right)\\
                                                                   &=& c^{n-1} \left( R_{\mu_{aa^{\ast}}} \freestar Moeb \right) + c^{n-1} \left( R_{\mu_{bb^{\ast}}} \freestar Moeb \right)\\
                                                                   &=& R_{\mu_{paa^{\ast}p}}^{p{\cal A}p} \freestar R_{\mu_c}^{-1} + R_{\mu_{pbb^{\ast}p}}^{p{\cal A}p} \freestar R_{\mu_c}^{-1}
\end{array}
\]
Using (\ref{boxedconvusefulform}), this can be written
\[
  R_{ \mu_{p(a+b)(a+b)^{\ast}p} \boxslash \mu_c } = R_{ \mu_{paa^{\ast}p} \boxslash \mu_c } + R_{ \mu_{pbb^{\ast}p} \boxslash \mu_c },
\]
which can equivalently be stated as
\[
  \left( \mu_{p(a+b)(a+b)^{\ast}p} \boxslash \mu_c \right) = \left( \mu_{paa^{\ast}p} \boxslash \mu_c \right) \boxplus \left( \mu_{pbb^{\ast}p} \boxslash \mu_c \right),
\]
which is what we had to prove.

Note that if $c\neq d$, $R_{\mu_{xx^{\ast}}} \freestar Moeb$ does
not appear as a factor in (\ref{cequalsd}). Therefore, there is no
reason why the result should hold for other Mar\u{c}henko Pastur
laws than $\mu_c$, since (\ref{connection2}) can not be used in
such cases. 

Although the proof of theorem~\ref{teo1d} is described for $c\leq 1$ only, 
the methods used in the proof of theorem~\ref{prop2} here can help us prove the case for $c>1$ also. 
If $R_n$ and $X_n$ are the random matrices from theorem~\ref{teo1d} and $c>1$, note that 
\[
  M_{\mu_{ \frac{1}{N}R_nR_n^{\ast} }} = 
  \frac{1}{c} M_{\mu_{ \frac{1}{N}R_n^{\ast}R_n }} = 
  \frac{1}{c} M_{\mu_{ \frac{1}{c} \frac{1}{n}R_n^{\ast}R_n }} =
  c^{-m-1} M_{\mu_{ \frac{1}{n}R_n^{\ast}R_n }},
\]
where $c^mf$ denoted the power series defined by $[ coef_k] (c^mf) = c^k [ coef_k](f)$. 
From this one can show that 
\begin{equation} \label{showequal}
\begin{array}{lll}
  R_{ \mu_{ \frac{1}{N}R_nR_n^{\ast} } } \freestar R_{\mu_c}^{-1} &=& M_{\mu_{ \frac{1}{N}R_nR_n^{\ast} }} \freestar Moeb \freestar (c^{m-1} Moeb) \\
                                                                  &=& \left( c^{-m-1} M_{\mu_{ \frac{1}{n}R_n^{\ast}R_n }} \right) \freestar Moeb \freestar (c^{m-1} Moeb)\\
                                                                  &=& \left( c^{-m-1} \left(  (c^{-m}) M_{\mu_{ \frac{1}{n}R_n^{\ast}R_n }}) \freestar (c^m Moeb) \right) \right) \freestar Moeb\\
                                                                  &=& \left( c^{-m-1} (M_{\mu_{ \frac{1}{n}R_n^{\ast}R_n }} \freestar Moeb) \right) \freestar \left( c^{-1} (c^{-m+1} Moeb) \right)\\
                                                                  &=& c^{-m} c^{-m-1} \left( M_{\mu_{ \frac{1}{n}R_n^{\ast}R_n }} \freestar Moeb \freestar (c^{-m+1}Moeb) \right)\\
                                                                  &=& c^{-2m-1} \left( R_{\mu_{ \frac{1}{n}R_n^{\ast}R_n }} \freestar R_{\mu_{\frac{1}{c}}}^{-1} \right).
\end{array}                                                                  
\end{equation}
Since theorem~\ref{teo1d} has been proved for $c\leq 1$, 
$R_{\mu_{ \frac{1}{n}(R_n+X_n)^{\ast}(R_n+X_n) }} \freestar R_{\mu_{\frac{1}{c}}}^{-1}$ will converge to 
\[
  R_{\mu_{ \frac{1}{n}R_n^{\ast}R_n }} \freestar R_{\mu_{\frac{1}{c}}}^{-1} + R_{\mu_{ \frac{1}{n}X_n^{\ast}X_n }} \freestar R_{\mu_{\frac{1}{c}}}^{-1}
\]
as $n\rightarrow\infty$ ($c>1$), and the result for $c>1$ follows from (\ref{showequal}). 
\section{Equivalence with known expressions for limit distributions of Information-Plus-Noise Type Matrices} \label{sectionsystem1}
~\citep{paper:raoedelman} studies systems where the sample covariance matrix is formed by taking independent samples of
a system of the form
\[
  y_n = A_nx_n + \sigma w_n
\]
(the noise factor $\sigma$ does not appear
in~\citep{paper:raoedelman}) where $x_n$ and $w_n$ are independent
standard (zero mean, unit variance) Gaussian random vectors, and $A_n$ is an $n\times L$ matrix.
The covariance matrix of the system is
$\Theta_n=A_nA_n^{\ast}+\sigma^2I$. In particular, when there is no
noise (i.e. $\sigma=0$), the covariance is $\Theta_n=A_nA_n^{\ast}$.
Denote by $\mu_{\Theta}$ the limiting eigenvalue distribution of
$A_nA_n^{\ast}$. \citep{paper:raoedelman} states that the limiting
eigenvalue distribution of the sample covariance matrix of the
system is $\left( \mu_{\Theta}\boxplus\mu_{\sigma^2I} \right)
\boxtimes \mu_c$. When there is no noise, the limit is
$\mu_{\Theta} \boxtimes \mu_c$. 
This way of passing from
\[
  \mu_{\Theta} \boxtimes \mu_c \mbox{ to } \left( \mu_{\Theta}\boxplus\mu_{\sigma^2I} \right) \boxtimes \mu_c
\]
is of course compatible with theorem~\ref{teo1d}. 
We will also show that it is equivalent with the results in~\citep{paper:doziersilverstein1}.
The following restrictions taken from~\citep{paper:doziersilverstein1} will be used:
\begin{enumerate}
  \item For $n=1,2,\cdots ,$, $X_n=(X_{ij}^n$, $n\times N$, i.d. for all $i,j,n$, independent across $i,j$ for each $n$, and $E|X_{11}^1-EX_{11}^1|^2 = 1$
  \item $R_n$ is $n\times N$ and independent of $X_n$, with $F^{\mu_{\Gamma_n}} \stackrel{\cal D}{\rightarrow} F^{\mu_{\Gamma}}$
\end{enumerate}

Theorem 1.1 in~\citep{paper:doziersilverstein1} expresses a relationship for
finding the limiting eigenvalue distribution $\mu_W$ of $W_n$ from that
of $\Gamma_n=\frac{1}{N}R_nR_n^{\ast}$ (denoted $\mu_{\Gamma}$).
More precisely, under the conditions 1) and 2), we have that (in a slightly rewritten form)
$F^{\mu_{W_n}} \stackrel{\cal D}{\rightarrow} F^{\mu_W}$ almost surely, where $F^{\mu_W}$ is a nonrandom p.d.f. characterized by
\begin{equation} \label{dozsilv11}
  m_{\mu_W}(z) = \int \frac{dF^{\mu_{\Gamma}}(t)}{ \frac{t}{1+\sigma^2 c m_{\mu_W}(z)} - (1+\sigma^2 c m_{\mu_W}(z))z + \sigma^2(1-c)}
\end{equation}
for any $z\in \C^+$.
The connection with multiplicative free convolution is hard to see from this formula.
To obtain this connection, the following lemma is needed, which will be proved in section~\ref{appendixa1}:
\begin{lem} \label{lem1}
(\ref{dozsilv11}) is equivalent to
\begin{equation} \label{doziersilv}
  m_{\mu_W}^{-1}\left( \frac{z}{1-\sigma^2 cz} \right) = (1-\sigma^2 cz)^2 m_{\mu_{\Gamma}}^{-1}(z) + \sigma^2(1-c) (1-\sigma^2 cz),
\end{equation}
for $z$ in some interval $(0,z_1)$.
\end{lem}
In (\ref{doziersilv}) and all other places where the inverse of the Stieltjes transform is taken in this paper,
we will mean the unique inverse on the negative real line. The inverse will only be calculated for positive values close to 0.
It will turn out that (\ref{doziersilv}) can be more conveniently expressed in terms of
distributions obtained from multiplicative free deconvolution with the Mar\u{c}henko Pastur law using the following lemma, 
which will be proved in section~\ref{appendixa2}:
\begin{lem} \label{lem2}
If
\begin{equation} \label{freedeconv}
  \mu_{\Gamma} = \mu_{\Theta} \boxtimes \mu_c,
\end{equation}
then, for $z$ in some interval $(0,z_1)$,
\begin{equation} \label{etatransforms}
  \eta_{\mu_{\Gamma}}^{-1}(z) = \frac{\eta_{\mu_{\Theta}}^{-1}(z)}{1-c+cz}
\end{equation}
and also
\begin{equation} \label{deconvolved}
  m_{\mu_{\Gamma}}^{-1}\left( \frac{z}{1-c-czm_{\mu_{\Theta}}^{-1}(z)} \right) = m_{\mu_{\Theta}}^{-1}(z) (1-c-czm_{\mu_{\Theta}}^{-1}(z)).
\end{equation}
\end{lem}

Using (\ref{deconvolved}), the following relationship with multiplicative free
convolution will be shown:
\begin{teo} \label{teo1a}
  Under the conditions 1) and 2), assume that
  \begin{equation}
    F^{\mu_{\Gamma_n}} \stackrel{\cal D}{\rightarrow} F^{\mu_{\Theta} \boxtimes \mu_c} \mbox{ a.s.}
  \end{equation}
  Then
  \begin{equation}
    F^{\mu_{W_n}} \stackrel{\cal D}{\rightarrow} F^{\left( \mu_{\Theta}\boxplus\mu_{\sigma^2I} \right) \boxtimes \mu_c} \mbox{ a.s}.
  \end{equation}
  Equivalently, assume that
  \begin{equation}
    F^{\mu_{\Gamma_n}} \stackrel{\cal D}{\rightarrow} F^{\mu_{\Gamma}} \mbox{ a.s.}
  \end{equation}
  Then
  \begin{equation}
    F^{\mu_{W_n}} \stackrel{\cal D}{\rightarrow} F^{\mu_W} \mbox{ a.s.}
  \end{equation}
  where $\mu_W$ is uniquely identified by the equation
  \begin{equation}
    \mu_W \boxslash \mu_c = (\mu_{\Gamma} \boxslash \mu_c) \boxplus \mu_{\sigma^2 I},
  \end{equation}
\end{teo}
Theorem~\ref{teo1a} will be proved in section~\ref{appendixa3}.

\subsection{The proof of lemma~\ref{lem1}} \label{appendixa1}
Rewritten in terms of the Stieltjes transform, (\ref{dozsilv11}) says that (with terms somewhat regrouped)
\[
  \frac{m_{\mu_W}}{1+\sigma^2 c m_{\mu_W}} = m_{\mu_{\Gamma}}\left( (1+\sigma^2 c m_{\mu_W})^2 z - \sigma^2 (1-c) (1+\sigma^2 c m_{\mu_W})\right) ,
\]
where $m_{\mu_W}$, $m_{\mu_{\Gamma}}$ are evaluated in $z$ when the parameter is omitted. 
We restrict ourselves to $z$ on the negative real line. 
The relation holds also for such $z$, since we can analytically continue to the negative real line. 
Evaluating in $m_{\mu_W}^{-1}(z)$ we get
\[
  \frac{z}{1+\sigma^2 c z} = m_{\mu_{\Gamma}}\left( (1+\sigma^2 c z)^2 m_{\mu_W}^{-1}(z) - \sigma^2 (1-c) (1+\sigma^2 c z)\right)
\]
for $z$ in some interval $(0,z_1)$. 
We will find it convenient to work with the inverse of the Stieltjes transform, so we rewrite the expression to 
\[
  m_{\mu_{\Gamma}}^{-1}\left( \frac{z}{1+\sigma^2 c z} \right) = (1+\sigma^2 c z)^2 m_{\mu_W}^{-1}(z) - \sigma^2 (1-c) (1+\sigma^2 c z) .
\]
Subsituting $u=\frac{z}{1+\sigma^2 cz}$ (or equivalently $z=\frac{u}{1-\sigma^2 cu}$) 
(this is an isomorphism of the positive real axis which sends $0$ to $0$), we get
\[
  m_{\mu_{\Gamma}}^{-1}(z) = \frac{ m_{\mu_W}^{-1}\left( \frac{z}{1-\sigma^2 cz} \right) }{(1-\sigma^2 cz)^2} - \frac{\sigma^2(1-c)}{1-\sigma^2 cz}
\]
for $z$ in some interval $(0,z_1)$, so that
\[
  m_{\mu_W}^{-1}\left( \frac{z}{1-\sigma^2 cz} \right) = (1-\sigma^2 cz)^2 m_{\mu_{\Gamma}}^{-1}(z) + \sigma^2(1-c) (1-\sigma^2 cz),
\]
which is (\ref{doziersilv}).

\subsection{The proof of lemma~\ref{lem2}} \label{appendixa2}
By the multiplicative property of the S-transform we have 
\[ 
  S_{\mu_{\Gamma}}(z) = \frac{S_{\mu_{\Theta}}(z)}{1+cz}. 
\]
Expressed in terms of the $\eta$-transform this can be written 
\[
  \eta_{\mu_{\Gamma}}^{-1}(z) = \frac{\eta_{\mu_{\Theta}}^{-1}(z)}{1-c+cz}, 
\]
which is (\ref{etatransforms}). 
Evaluating in $\eta_{\mu_{\Theta}}(z)$ and applying $\eta_{\mu_{\Gamma}}$ on both sides gives
\[
  \eta_{\mu_{\Gamma}}\left( \frac{z}{1-c+c\eta_{\mu_{\Theta}}(z)} \right) = \eta_{\mu_{\Theta}}(z)
\]
for $z\geq 0$. 
This can also be expressed in terms of Stieltjes transforms as
\[
  \frac{m_{\mu_{\Gamma}}\left( -\frac{1-c+c\eta_{\mu_{\Theta}}(z)}{z} \right)}{\frac{z}{1-c+c\eta_{\mu_{\Theta}}(z)}} = \frac{m_{\mu_{\Theta}}(-\frac{1}{z})}{z}
\]
Regrouping terms and substituting $-\frac{1}{z}$ for $z$ we get
\[
  m_{\mu_{\Gamma}}\left( z(1-c-cz m_{\mu_{\Theta}}(z)) \right) = \frac{m_{\mu_{\Theta}}(z)}{1-c-czm_{\mu_{\Theta}}(z)}
\]
for $z<0$. 
Substituting $m_{\mu_{\Theta}}^{-1}(z)$ for $z$ and taking the inverse Stieltjes transform $m_{\mu_{\Gamma}}^{-1}$ we get
\[
  m_{\mu_{\Gamma}}^{-1}\left( \frac{z}{1-c-czm_{\mu_{\Theta}}^{-1}(z)} \right) = m_{\mu_{\Theta}}^{-1}(z) (1-c-czm_{\mu_{\Theta}}^{-1}(z))
\]
for $z$ in some interval $(0,z_1)$, 
which is (\ref{deconvolved}).

\subsection{The proof of theorem~\ref{teo1a}} \label{appendixa3}
Note that if 
\[
  z = \frac{z_1}{1-c-cz_1 m_{\mu_{\Theta}}^{-1}(z_1)} 
\]
for $z_1$ positive and close to $0$, 
then we have 
\begin{equation} \label{substitute}
  \frac{z}{1-\sigma^2 cz} = \frac{z_1}{1-c-cz_1 \left( m_{\mu_{\Theta}}^{-1}(z_1) + \sigma^2 \right)}.
\end{equation}
Note also that $\frac{z}{1-\sigma^2 cz}\geq 0$ as long as $z < \frac{1}{\sigma^2 c}$. 
Substituting (\ref{substitute}) and (\ref{deconvolved}) in (\ref{doziersilv}) we get $m_{\mu_W}^{-1}\left( \frac{z}{1-\sigma^2 cz} \right)=$
\begin{equation} \label{calculation}
  \begin{array}{lll}
    &   & \left( 1 - \frac{\sigma^2 cz_1}{1-c-cz_1 m_{\mu_{\Theta}}^{-1}(z_1)} \right)^2 m_{\mu_{\Theta}}^{-1}(z_1) (1-c-cz_1 m_{\mu_{\Theta}}^{-1}(z_1)) \\
    &   & + \sigma^2(1-c)  \left( 1 - \frac{\sigma^2 cz_1}{1-c-cz_1 m_{\mu_{\Theta}}^{-1}(z_1)} \right) \\
    & = & \frac{ \left( 1-c-cz_1 m_{\mu_{\Theta}}^{-1}(z_1)-\sigma^2 cz_1 \right)^2 m_{\mu_{\Theta}}^{-1}(z_1) + \sigma^2(1-c) \left( 1-c-cz_1 m_{\mu_{\Theta}}^{-1}(z_1)-\sigma^2 cz_1 \right) }{1-c-cz_1 m_{\mu_{\Theta}}^{-1}(z_1)}\\
    & = & \frac{ \left( 1-c-cz_1 m_{\mu_{\Theta}}^{-1}(z_1)-\sigma^2 cz_1 \right) \left( \left( 1-c-cz_1 m_{\mu_{\Theta}}^{-1}(z_1)-\sigma^2 cz_1 \right) m_{\mu_{\Theta}}^{-1}(z_1) + \sigma^2(1-c) \right) }{1-c-cz_1 m_{\mu_{\Theta}}^{-1}(z_1)}\\
    & = & \frac{ \left( 1-c-cz_1 (m_{\mu_{\Theta}}^{-1}(z_1)+\sigma^2) \right) \left( 1-c-cz_1 m_{\mu_{\Theta}}^{-1}(z_1) \right) \left( m_{\mu_{\Theta}}^{-1}(z_1)+\sigma^2 \right)}{1-c-cz_1 m_{\mu_{\Theta}}^{-1}(z_1)}\\
    & = & \left( m_{\mu_{\Theta}}^{-1}(z_1)+\sigma^2 \right) \left( 1-c-cz_1 (m_{\mu_{\Theta}}^{-1}(z_1)+\sigma^2) \right) \\
    & = & m_{\mu_{\Theta}\boxplus\mu_{\sigma^2 I}}^{-1}(z_1) \left( 1-c-cz_1 m_{\mu_{\Theta}\boxplus\mu_{\sigma^2 I}}^{-1}(z_1) \right) 
  \end{array}
\end{equation}
Here we have used that $m_{\mu_{\Theta}\boxplus\mu_{\sigma^2 I}}^{-1}(z_1) = m_{\mu_{\Theta}}^{-1}(z_1) +\sigma^2$, 
which follows from the additivity property of the $R$-transform and 
the fact that the inverse of the Stieltjes transform is used to define the $R$-transform 
(perform additive free convolution with $\mu_{\sigma^2I}$). 
(\ref{calculation}) is thus nothing else than (\ref{deconvolved}) (with $m_{\mu_{\Theta}}$ replaced by $m_{\mu_{\Theta}\boxplus\mu_{\sigma^2 I}}$). 
Since (\ref{deconvolved}) is just an equivalent expression for multiplicative free convolution, we therefore have 
\[
  \mu_W = \left( \mu_{\Theta}\boxplus\mu_{\sigma^2 I} \right) \boxtimes \mu_c,
\]
or equivalently 
\[
  \mu_W \boxslash \mu_c = \mu_{\Theta}\boxplus\mu_{\sigma^2 I} = \left( \mu_{\Gamma} \boxslash \mu_c \right) \boxplus \mu_{\sigma^2 I}.
\]
This completes the proof.

\section{Using $G$-analysis to estimate the spectral function of covariance matrices} \label{ganalysis}
It turns out that multiplicative free deconvolution can also be
used to estimate covariance matrices. The general statistical
analysis of observations, also called {\em
$G$-analysis}~\citep{book:girkostat} is a mathematical theory for
complex systems where the number of parameters of the underlying
mathematical model increase together with the growth of the number
of observations of the system. The mathematical models which
approach the system in some sense are called {\em $G$-estimators}.
The main difficulty in $G$-analysis is to find good
$G$-estimators. $G$-estimators have already shown their usefulness
in many applications~\citep{paper:mestre}. We denote by $N$ the
number of observations of the system, and by $n$ the number of
parameters of the mathematical model. The condition used in
$G$-analysis expressing the growth of the number of observations
vs. the number of parameters in the mathematical model, is called
the {\em $G$-condition}. The $G$-condition used throughout this
paper is (\ref{gcondition}).

Girko restricts to systems where a number of independent random
vector observations are taken, and where the random vectors have
identical distributions. If a random vector $r_n$ has length $n$,
we will let $\Theta_n$ denote it's covariance, while $\Gamma_n$
will still denote sample covariance matrices. The $\Gamma_n$ we
analyze in this section are more restrictive than in previous
sections, since independence across samples is assumed. Girko
calls estimators for the Stieltjes transform of covariance
matrices {\em $G^2$-estimators}. In chapter 2.1
of~\citep{chapter:girkotenyears} he introduces the following
expression as candidate for a $G^2$-estimator:
\begin{equation} \label{gestcondition0}
  G^2_n(z) = \frac{\hat{\theta}(z)}{z} m_{\mu_{\Gamma_n}}(\hat{\theta}(z)),
\end{equation}
where the function $\hat{\theta}(z)$ is the solution to the equation
\begin{equation} \label{gestcondition}
  \hat{\theta}(z) c m_{\mu_{\Gamma_n}}(\hat{\theta}(z)) - (1-c) + \frac{\hat{\theta}(z)}{z} = 0.
\end{equation}
Girko claims that a function $G^2_n(z)$ satisfying (\ref{gestcondition}) and (\ref{gestcondition0}) is a good approximation for
the Stieltjes transform of the covariance matrices
$m_{\Theta_n}(z) = tr_n \left\{ \Theta_n - z I_n \right\}^{-1}$.
More precisely, he shows that when (\ref{gestcondition0}), (\ref{gestcondition}) and the $G$-condition (\ref{gcondition}) are fulfilled,
under certain conditions there exists a $c>0$ such that
\begin{equation} \label{limvals}
  \lim_{n\rightarrow\infty} \sup_{\stackrel{0 < c \leq \Im(z) \leq S}{| \Re(z) | \leq T}} \left| G^2_n(z) - m_{\Theta_n}(z) \right| = 0,
\end{equation}
with probability one for every $S>0$ and $T>0$.
According to Girko, analytical continuation of $G^2_n(z)$ can be performed to
obtain limits for other $z$ than the ones in (\ref{limvals}).

As it turns out, the $G^2$-estimator can equivalently be expressed in terms of multiplicative free convolution:
\begin{teo} \label{teo2}
  For the $G^2$-estimator given by (\ref{gestcondition0}), (\ref{gestcondition}), the following holds for real $z < 0$:
  \begin{equation}
    G^2_n(z) = m_{\mu_{\Gamma_n} \boxslash \mu_c}
  \end{equation}
\end{teo}

\begin{proof}
(\ref{gestcondition}) can be rewritten to
\[
  -c\eta_{\mu_{R_n}} \left( -\frac{1}{\hat{\theta}(z)} \right) - (1-c) + \frac{\hat{\theta}(z)}{z} = 0
\]
\[
  \eta_{\mu_{R_n}}^{-1}\left( \frac{1}{c}\left( \frac{\hat{\theta}(z)}{z} -(1-c) \right) \right) = -\frac{1}{\hat{\theta}(z)},
\]
which we will write
\begin{equation} \label{gestcondition2}
  - \frac{1}{\eta_{\mu_{R_n}}^{-1}\left( \frac{1}{c}\left( \frac{\hat{\theta}(z)}{z} -(1-c) \right) \right)} = \hat{\theta}(z).
\end{equation}
Denote by $\mu$ the measure with Stieltjes transform $G^2_n(z)$.
(\ref{gestcondition0}) can be rewritten using the $\eta$-transform as
\[
  \eta_{\mu}\left( -\frac{1}{z} \right) = \eta_{\mu_{R_n}}\left( -\frac{1}{\hat{\theta}(z)} \right).
\]
Since $\eta_{\mu_{R_n}}$ and $\eta_{\mu}$ are monotone, it is easily seen from this that $\hat{\theta}$ is monotone since it is a combination of
monotone functions.
Forming the inverse functions on both sides, and also applying $\hat{\theta}$, yields
\begin{equation} \label{etaconnection}
  \hat{\theta}\left( -\frac{1}{\eta_{\mu}^{-1}(z)} \right) = -\frac{1}{\eta_{\mu_{R_n}}^{-1}(z)}.
\end{equation}
Showing $\mu_{R_n} = \mu \boxtimes \mu_c$ is equivalent to (after rearranging (\ref{etatransforms}))
\[
  -\frac{1}{\eta_{\mu}^{-1}(z)} = -\frac{1}{(1-c+cz)\eta_{\mu_{R_n}}^{-1}(z)}
\]
Applying $\hat{\theta}$ on both sides and using (\ref{etaconnection}) yields that this is equivalent to
\begin{equation} \label{vshs}
  -\frac{1}{\eta_{\mu_{R_n}}^{-1}(z)} = \hat{\theta}\left( -\frac{1}{(1-c+cz)\eta_{\mu_{R_n}}^{-1}(z)} \right)
\end{equation}

Observe now that (\ref{vshs}) and (\ref{gestcondition2}) are related in the following way:
If we substitute $z = \frac{1}{c}\left( \frac{\hat{\theta}(w)}{w} - (1-c)\right)$ into (\ref{vshs}), the argument on the right hand side can be rewritten using (\ref{gestcondition2}) to
\[
    -\frac{1}{
               \left(
                 1-c+\frac{\hat{\theta}(z)}{z} -(1-c)
               \right)
               \left(
                 -\frac{1}{\hat{\theta}(z)}
               \right)
             }
    = z,
\]
so that (\ref{vshs}) is nothing but a restatement of (\ref{gestcondition2}),
at least on values of the form $z = \frac{1}{c}\left( \frac{\hat{\theta}(w)}{w} - (1-c)\right)$.
If these values take on an open set of real values,
equality in (\ref{vshs}) follows for all $z$ by analytic continuation.
This happens when $\hat{\theta}(z)\neq kz$ for some constant $k$.
If $\hat{\theta}(z) = kz$, then $\eta_{\mu_{R_n}}$ is seen to be constant, which only happens in trivial cases.
Thus we have that $\mu_{R_n} = \mu \boxtimes \mu_c$, and we are done
\end{proof}

Several remarks concerning theorem~\ref{teo2} are in place.
First of all, the $G^2$-estimator has a much shorter expression in terms of multiplicative free deconvolution, which also 
places it as an ingredient in theorem~\ref{teo1d}. 
The theorem is nice to combine with continuity results for free convolution.
Voiculescu has proved such results when convergence is in the weak-$\ast$ topology~\citep{paper:vounbounded}.
This enables us in many cases to conclude that
\[
  \lim_{n\rightarrow\infty} G^2_n(z) = \lim_{n\rightarrow\infty} m_{\mu_{\Gamma_n} \boxslash \mu_c} = m_{\mu_{\Gamma} \boxslash \mu_c}
\]
for some probability measure $\mu_{\Gamma}$.
Secondly,~\citep{paper:raoedelman} expresses the exact same estimator, i.e.
\[
  \lim_{n\rightarrow\infty} \mu_{\Gamma_n} \boxslash \mu_c = \lim_{n\rightarrow\infty} \mu_{\Theta_n}
\]
in the case of Gaussian systems. Theorem~\ref{teo2} can be seen as a way of generalizing from the Gaussian case.

\section{Further work}
The concept of freeness and free convolution can be extended to
unbounded random variables and general probability measures.
In~\citep{paper:vounbounded} it is shown how this can be done in
the context of unbounded operator spaces, and certain regularity
properties are proved. For instance, if $\mu_n\rightarrow\mu$ and
$\nu_n\rightarrow\nu$ in the weak-$\ast$ topology with both
$\mu\neq\delta_0$ and $\nu\neq\delta_0$, then
$\mu_n\boxtimes\nu_n\rightarrow\mu\boxtimes\nu$ in the weak-$\ast$
topology also. It is possible that applying such extensions
together with the methods applied here can extend the results to the 
same generality as those in~\citep{paper:doziersilverstein1}. 
This may be addressed in a future paper. 

The $G^2$-estimator is just one of many estimators introduced by Girko. 
He has estimators for many other quantities also~\citep{chapter:girkotenyears}, like for the square root and the moments of covariance matrices. 
Certain of these estimators may also have alternative expressions in terms of free probability constructs. 

\bibliography{../bib/mybib,../bib/mainbib,../bib/drbib}
\end{document}